\documentclass[journal]{IEEEtran}

\usepackage{graphicx}
\usepackage{color}
\usepackage{amsfonts}
\usepackage{float}
\usepackage{amsmath}
\usepackage{amsthm}
\usepackage{epstopdf}
\usepackage{changepage}
\usepackage{latexsym}
\usepackage{tikz}
\usepackage{cases}
\usepackage{longtable}
\usepackage{supertabular}
\usepackage{tabu}
\usepackage{multirow}
\usepackage{multicol}
\usepackage{booktabs}
\usepackage{arydshln}
\usepackage{setspace}
\usepackage{algpseudocode}
\usepackage{etex}
\usepackage{cite}
\usepackage{balance}

\usepackage{cancel}

\IEEEoverridecommandlockouts
\usepackage{graphics} 
\usepackage{epsfig} 
\usepackage{times} 

\usepackage{algorithm}

\makeatletter
\def\BState{\State\hskip-\ALG@thistlm}
\makeatother

\usepackage{amssymb}

\usepackage{amsthm}
\usepackage{amsmath} 
\usepackage{amssymb}  
\usepackage{epstopdf}
\usepackage{epsfig}
\usepackage{subfigure}
\newcommand{\m}{\boldsymbol}
\usepackage{tikz}
\usetikzlibrary{dsp,fit}
\usetikzlibrary{decorations.pathmorphing} 

\usetikzlibrary{dsp,fit}
\usetikzlibrary{decorations.pathmorphing} 
\tikzstyle{block} = [draw,rectangle,thick,minimum height=3em,minimum width=5em]
\tikzstyle{branch} = [circle,inner sep=1pt,minimum size=2mm,fill=white,draw=black]
\tikzstyle{branch1} = [inner sep=1pt,minimum size=0.1mm,fill=white,draw=black]
\tikzstyle{branch1} = [circle,inner sep=0pt,minimum size=1mm,fill=black,draw=black]

\makeatletter
\dspdeclareoperator{dspvoidshapeadder}{
    \pgfutil@tempdima=\radius
    \pgfutil@tempdima=0.55\pgfutil@tempdima
    \pgfusepathqstroke
}
\tikzstyle{connector2} = [-,thick]
\tikzstyle{connector} = [->,thick]
\tikzstyle{connector1} = [->,ultra thick]
\tikzstyle{snakeline} = [connector, decorate, decoration={pre length=0.2cm,
                         post length=0.2cm, snake, amplitude=.7mm,
                         segment length=2mm},very thick, red, ->]
                         \tikzstyle{snakeline1} = [connector, decorate, decoration={pre length=0.2cm,
                                                  post length=0.2cm, snake, amplitude=.3mm,
                                                  segment length=2mm},very thick, blue, ->]
  \tikzstyle{snakeline2} = [connector, decorate, decoration={pre length=0.4cm,
                                                  post length=0.4cm, snake, amplitude=1mm,
                                                  segment length=4mm},very thick, red, ->]
\tikzset{
vdspadder/.style={
  shape=dspvoidshapeadder,
   line cap=rect,
  line join=rect,
  line width=\dspblocklinewidth,
  minimum size=\dspoperatordiameter,
  label={185:$+$},
  label={265:$-$}
  },
vadspadder/.style={
  shape=dspvoidshapeadder,
  line cap=rect,
  line join=rect,
  line width=\dspblocklinewidth,
  minimum size=\dspoperatordiameter,
  label=below right:$-$,
  label=above right:$+$
  }
}
\makeatother
\usepackage{amsmath}
\usepackage{amssymb}
\usepackage{textcase}
\usepackage{soul}
\usepackage{indentfirst}

\newcommand{\normof}[1]{\|#1\|}

\newcommand{\bmat}[1]{\begin{bmatrix} #1 \end{bmatrix}}

\newcommand*{\LargerCdot}{\raisebox{-0.25ex}{\scalebox{1.5}{$\cdot$}}}
\newcommand{\tabincell}[2]{\begin{tabular}{@{}#1@{}}#2\end{tabular}}

\begin{document}

\title{Dynamic State Estimation for Multi-Machine Power System by Unscented Kalman Filter with Enhanced Numerical Stability}

\author{Junjian~Qi,~\IEEEmembership{Member,~IEEE, }
        Kai~Sun,~\IEEEmembership{Senior Member,~IEEE, }%
        Jianhui~Wang,~\IEEEmembership{Senior Member,~IEEE,}
        and Hui Liu,~\IEEEmembership{Member,~IEEE}
        \thanks{This work was supported in part by U.S. Department of Energy, Office of Electricity Delivery and Energy Reliability and the CURENT Engineering Research Center. 
        Paper no. TSG-00759-2015.
        
        J.~Qi and J. Wang are with the Energy Systems Division, Argonne National Laboratory, Argonne, IL 60439 USA (e-mails: jqi@anl.gov; jianhui.wang@anl.gov).

        K. Sun is with the Department of Electrical Engineering and Computer Science, University of Tennessee, Knoxville, TN 37996 USA (e-mail: kaisun@utk.edu).
        
        H. Liu is with the Department of Electrical Engineering, Guangxi University, Nanning, 530004 China and was a visiting scholar at the Energy Systems Division, 
        Argonne National Laboratory, Argonne, IL 60439 USA (e-mail: hughlh@126.com).
        
}
}

\markboth{preprint of doi: 10.1109/TSG.2016.2580584, IEEE Transactions on Smart Grid}{stuff}
\maketitle

\begin{abstract}
In this paper, in order to enhance the numerical stability of the unscented Kalman filter (UKF) used for power system dynamic state estimation, 
a new UKF with guaranteed positive semidifinite estimation error covariance (UKF-GPS) is proposed and  
compared with five existing approaches, including UKF-schol, UKF-$\kappa$, UKF-modified, UKF-$\Delta Q$, and the square-root unscented Kalman filter (SR-UKF). 
These methods and the extended Kalman filter (EKF) are tested by performing dynamic state estimation on WSCC 3-machine 9-bus system and NPCC 48-machine 140-bus system. 
For WSCC system, all methods obtain good estimates. However, for NPCC system, both EKF and the classic UKF fail. 
It is found that UKF-schol, UKF-$\kappa$, and UKF-$\Delta Q$ do not work well in some estimations while UKF-GPS works well in most cases. 
UKF-modified and SR-UKF can always work well, indicating their better scalability mainly due to the enhanced numerical stability.
\end{abstract}


\section*{Nomenclature}
\addcontentsline{toc}{section}{Nomenclature}
\begin{IEEEdescription}[\IEEEusemathlabelsep\IEEEsetlabelwidth{$V_1,V_2,V_3,V_4$}]
\item [$\boldsymbol{0}_{a,b}$] Zero matrix with dimension $a\times b$.
\item [$\boldsymbol{I}_a$]  Identity matrix with dimension $a$. 
\item [$\boldsymbol{f}_c, \boldsymbol{f}$]  Column vector of continuous and discrete state transition functions.
\item [$\boldsymbol{h}_c, \boldsymbol{h}$] Column vector of continuous and discrete measurement functions.
\item [$\boldsymbol{K}$] Kalman gain matrix. 
\item [$\boldsymbol{m}$] Estimated mean of the state. 
\item [$\boldsymbol{m}_0, \boldsymbol{m}^-$] Initial and predictd mean of the state. 
\item [$\boldsymbol{P}_0, \boldsymbol{P}^-, \boldsymbol{P}$] Initial, predicted and updated estimation error covariance. 
\item [$\boldsymbol{P}_{\tilde{\boldsymbol{y}}_k\tilde{\boldsymbol{y}}_k}$] Covariance of the measurement. 
\item [$\boldsymbol{P}_{\boldsymbol{x}_k\boldsymbol{y}_k}$] Cross covariance of the state and measurement. 
\item [$\boldsymbol{q}, \boldsymbol{r}$] Process noise and measurement noise column vectors.
\item [$\boldsymbol{Q}, \boldsymbol{R}$] Constant covariance matrices of $\boldsymbol{q}$ and $\boldsymbol{r}$. 
\item [$\boldsymbol{S}$] Cholesky factor (matrix square root) of the estimation error covariance $\boldsymbol{P}$. 
\item [$\boldsymbol{W}_m, \boldsymbol{W}_c$] Weights for the mean and the covariance of the state or measurement. 
\item [$\boldsymbol{x}$] Column vector of the states. 
\item [$\boldsymbol{\mathcal{X}},\boldsymbol{\mathcal{X}}^-$] Sigma points and predicted sigma points. 
\item [$\boldsymbol{y}$] Column vector of the measurements. 
\item [$\boldsymbol{y}^-$] Predicted measurement. 
\item [$\boldsymbol{\mathcal{Y}}^-$] Propagated sigma points by the measurement function. 
\item [$\delta$] Rotor angle in rad. 
\item [$\omega, \omega_0$] Rotor speed and rated rotor speed in rad/s. 
\item [$\Psi$] Voltage source.
\item [$\boldsymbol{\Psi}_\textbf{R}, \boldsymbol{\Psi}_\textbf{I}$] Column vectors of all generators' real and imaginary parts of the voltage source on system reference frame. 
\item [$E_{\textrm{fd}}$] Internal field voltage in pu. 
\item [$E_\textrm{t}$] Terminal voltage phasor. 
\item [$e_\textrm{q}, e_\textrm{d}$] Terminal voltage at q and d axes in pu. 
\item [$e'_\textrm{q}, e'_\textrm{d}$] Transient voltage at q and d axes in pu. 
\item [$e_\textrm{R}, e_\textrm{I}$] Real and imaginary part of the terminal voltage phasor. 
\item [$e_x$] System state error averaged for one type of state ($\delta$, $\omega$, $e'_\textrm{q}$, or $e'_\textrm{d}$) over a time period. 
\item [$g$] Number of generators. 
\item [$g_2, g_4$] Number of generators with classical model and fourth-order transient model.
\item [$\bar{g}$] Number of PMUs. 
\item [$\mathcal{G}_2, \mathcal{G}_4$] Set of generators with second-order classical model and fourth-order transient model. 
\item [$\mathcal{G}_\textrm{P}$] Set of generators where PMUs are installed. 
\item [$H$] Generator inertia constant in second. 
\item [$I_\textrm{t}$] Terminal current phasor. 
\item [$i_\textrm{q}, i_\textrm{d}$] Current at q and d axes in pu.
\item [$i_\textrm{R}, i_\textrm{I}$] Real and imaginary part of the terminal current phasor in pu. 
\item [$K_\textrm{D}$] Damping factor in pu. 
\item [$n, v, p$] Number of states, inputs, and outputs. 
\item [$P_\textrm{e}$] Electrical active output power in pu. 
\item [$S_\textrm{B}, S_\textrm{N}$] System and generator base MVA. 
\item [$T_\textrm{m}, T_\textrm{e}$] Mechanical torque and electric air-gap torque in pu.
\item [$T'_{\textrm{q0}}, T'_{\textrm{d0}}$] Open-circuit time constants for q and d axes in second.
\item [$x_{\textrm{q}}, x_{\textrm{d}}$] Synchronous reactance at q and d axes in pu.
\item [$x'_{\textrm{q}}, x'_{\textrm{d}}$] Transient reactance at q and d axes in pu. 
\item [$\boldsymbol{\overline{Y}}$] Admittance matrix of the reduced network only consisting of generators\footnote{The elements of $\boldsymbol{\overline{Y}}$ are constant if the difference between $x'_\textrm{d}$ and $x'_\textrm{q}$ is ignored (\cite{Wang1}).}. 
\item [$\boldsymbol{\overline{Y}}_i$] The $i$th row of $\boldsymbol{\overline{Y}}$. 
\item [$\operatorname{schol}(\cdot)$] Cholesky factor of a matrix. 
\item [$\operatorname{cholupdate}(\cdot)$] Rank 1 update to Cholesky factorization. 
\item [$\operatorname{eig}(\cdot)$] Obtain the eigenvalue and eigenvector of a matrix. 
\item [$\operatorname{diag}(\cdot)$] Create diagonal matrix or get diagonal elements of matrix. 
\item [$\operatorname{qr}(\cdot)$] Orthogonal-triangular decomposition of a matrix. 
\item [$\sqrt{\boldsymbol{P}}$] Matrix square root of a positive semidefinite matrix $\boldsymbol{P}$, which is a matrix  $\boldsymbol{S}=\sqrt{\boldsymbol{P}}$ such that $\boldsymbol{P}=\boldsymbol{S}\boldsymbol{S}^\top$.
\item [$\operatorname{Re}(\cdot), \operatorname{Im}(\cdot)$] Real part and imaginary part. 
\item [${\left[\cdot\right]}_i$]  The $i$th column of a matrix.
\item [${\left[\cdot\right]}_\mathcal{S}$]  Columns of a matrix belonging to a set $\mathcal{S}$. 
\item [$||\cdot||$]   Frobenius norm of a matrix.
\item [$||\cdot||_2$]   2-norm of a vector.
\item [$\LargerCdot$, $\times$]  Elementwise product and matrix product.
\end{IEEEdescription}

\section{Introduction}

\IEEEPARstart{S}{tate} estimation is an important application of the energy management system (EMS). 
However, the widely studied static state estimation [\ref{se1}]--[\ref{se7}] assumes that the power system operates in quasi-steady state,
based on which the static states of the system, i.e. the voltage magnitude and phase angles of all buses, 
are estimated by making use of the Supervisory Control and Data Acquisition (SCADA) and/or phasor measurement unit (PMU) measurements. 

Static state estimation is important for power system monitoring and also provides input data for other important applications in EMS, 
such as automatic generation control (AGC) and optimal power flow (OPF). 
However, it may not be sufficient for good system monitoring and situational awareness 
as the power system becomes more dynamic due to the increasing penetration of renewable generation that has very high uncertainty and variation. 
Therefore, accurate dynamic states of the system obtained from real-time dynamic state estimation (DSE)
facilitated by high-level PMU deployment has thus become essential.
With the high global positioning system (GPS) synchronization accuracy, 
PMUs can provide highly synchronized measurements of voltage and current phasors in high sampling rate, 
thus playing a critical role in achieving real-time wide-area monitoring, protection, and control.

Power system DSE has been implemented by different types of Kalman filters. 
The most common application of the Kalman filter (KF) [\ref{kf}] to nonlinear systems is in the form of extended Kalman filter (EKF) [\ref{ekf1}], [\ref{ekf2}], 
which linearizes all nonlinear transformations and substitutes Jacobian matrices for the linear transformations in KF equations,  
based on the assumption that all transformations are quasi-linear. Power system DSE has been implemented by EKF [\ref{ekf3}], [\ref{ekf4}].

Although EKF maintains the elegant and computationally efficient recursive update form of the KF, it works well only in a `mild' nonlinear environment 
due to the first-order Taylor series approximation for nonlinear functions \cite{CKF}. 
It is sub-optimal and can easily lead to divergence. 
The linearized transformations are reliable only when the error propagation can be well approximated by a linear function. 
Also, the linearization can be applied only if the Jacobian matrix exists. 
Even if the Jocobian matrix esists, calculating it can be a difficult and error-prone. 

The unscented transformation (UT) [\ref{ut}] was developed to address the deficiencies of linearization by providing a more direct and explicit mechanism for transforming mean and covariance information. Based on UT, Julier et al. \cite{ukf,ukf2,ukf1} proposed the unscented Kalman filter (UKF) as a derivative-free alternative to EKF in the framework of state estimation. 
The UKF has been applied to power system DSE, for which no linearization or calculation of Jacobian matrices is needed \cite{pwukf1}, \cite{pwukf2}. 
However, in \cite{pwukf1} and \cite{pwukf2} UKF is only applied to estimate the dynamic states for the single-machine infinite-bus system or WSCC 3-machine system. 

It is not surprising that UKF has not been applied to larger power systems. As has been pointed out in \cite{CKF} and \cite{bellman}, 
both EKF and UKF can suffer from the curse of dimensionality and the effect of dimensionality may become detrimental in high-dimensional state-space models with state-vectors of size twenty or more, especially when there are high degree of nonlinearities in the equations that describe the state-space model, which is exactly the case for power systems. 

Therefore, even if classic UKF has good performance for small systems, it might not work at all for large power systems. 
We will show that it is the numerical stability that mainly limits the scalability of the classic UKF. 
Specifically, when the estimation error covariance is propagated, it sometimes cannot maintain the positive semidefiniteness, 
thus making its square-root unable to be calculated. 

In this paper, we introduce and compare six techniques that can be used to enhance the numerical stability of UKF, 
including the EKF/UKF toolbox approach \cite{ekfukf}, properly setting a parameter of unscented transformation \cite{sr_ukf1}, the modified UKF approach \cite{ukf2}, 
adding an extra positive definite matrix \cite{deltaQ,deltaQ1}, the UKF with guaranteed positive semidifinite estimation error covariance (UKF-GPS) proposed in this paper, and the square-root UKF (SR-UKF) \cite{sr_ukf}. 

The remainder of this paper is organized as follows. Section \ref{basis} briefly introduces the unscented transformation and the classic UKF procedure. 
Section \ref{ns} discusses six techniques for enhancing the numerical stability of the classic UKF. 
Section \ref{ps dse} explains how Kalman filters can be implemented for power system dynamic state estimation. 
Section \ref{result} tests the proposed methods on the WSCC 3-machine 9-bus system and NPCC 48-machine 140-bus system.
Finally the conclusion is drawn in Section \ref{conclusion}.

\section{Unscented Kalman Filter} \label{basis}

A discrete-time nonlinear system can be described as
\begin{subnumcases}{\label{non}}
\boldsymbol{x}_{k}=\boldsymbol{f}(\boldsymbol{x}_{k-1},\boldsymbol{u}_{k-1})+\boldsymbol{q}_{k-1} \\
\boldsymbol{y}_{k}=\boldsymbol{h}(\boldsymbol{x}_k,\boldsymbol{u}_k)+\boldsymbol{r}_k,
\end{subnumcases}
where $\boldsymbol{x}_k \in \mathbb{R}^n$, $\boldsymbol{u}_k \in \mathbb{R}^v$, and $\boldsymbol{y}_k \in \mathbb{R}^p$ are, respectively, 
state variables, inputs, and observed measurements at time step $k$; the estimated mean and estimation error covariance are $\boldsymbol{m}$ and $\boldsymbol{P}$; 
$\boldsymbol{f}$ and $\boldsymbol{h}$ are vectors consisting of nonlinear state transition functions and measurement functions; $\boldsymbol{q}_{k-1} \sim N(0,\boldsymbol{Q}_{k-1})$ is 
the Gaussian process noise at time step $k-1$; $\boldsymbol{r}_k \sim N(0,\boldsymbol{R}_k)$ is the Gaussian measurement noise at time step $k$; and 
$\boldsymbol{Q}_{k-1}$ and $\boldsymbol{R}_k$ are covariances of $\boldsymbol{q}_{k-1}$ and $\boldsymbol{r}_k$.

\subsection{Unscented Transformation}
Unscented Transformation (UT) is proposed based on the idea that ``it is easier to approximate a probability distribution than it is to approximate an arbitrary nonlinear function or transformation'' \cite{ut}. A set of sigma points are chosen so that their mean and covariance are $\boldsymbol{m}$ and $\boldsymbol{P}$. The nonlinear function is applied to each point to yield a cloud of transformed points and the statistics of the transformed points can then be calculated to form an estimate of the nonlinearly transformed mean and covariance. 

Specifically, a total of $2\,n+1$ sigma points (denoted by $\boldsymbol{\mathcal{X}}$) are calculated from the columns of the matrix $\eta\sqrt{\boldsymbol{P}}$ as 
\begin{subnumcases}{\label{sigma}}
\boldsymbol{\mathcal{X}}^{(0)}=\boldsymbol{m} \\
\boldsymbol{\mathcal{X}}^{(i)}=\boldsymbol{m}+\left[\eta\sqrt{\boldsymbol{P}}\right]_i,\quad i=1,\ldots,n   \\
\boldsymbol{\mathcal{X}}^{(i)}=\boldsymbol{m}-\left[\eta\sqrt{\boldsymbol{P}}\right]_i,\quad i=n+1,\ldots,2\,n 
\end{subnumcases}
with weights
\begin{subnumcases}{\label{weight}}
\boldsymbol{W}_m^{(0)}=\frac{\lambda}{n+\lambda}  \\
\boldsymbol{W}_c^{(0)}=\frac{\lambda}{n+\lambda}+(1-\alpha^2+\beta)  \\
\boldsymbol{W}_m^{(i)}=\frac{1}{2(n+\lambda)}, \quad i=1,\ldots,2\,n  \\
\boldsymbol{W}_c^{(i)}=\frac{1}{2(n+\lambda)}, \quad i=1,\ldots,2\,n,
\end{subnumcases}
where the matrix square root of a positive semidefinite matrix $\boldsymbol{P}$ is a matrix $\boldsymbol{S}=\sqrt{\boldsymbol{P}}$ such that $\boldsymbol{P}=\boldsymbol{S}\boldsymbol{S}^\top$, $\boldsymbol{W}_m$ and $\boldsymbol{W}_c$ are respectively weights for the mean and the covariance, $\eta=\sqrt{n+\lambda}$, $\lambda$ is a scaling parameter defined as $\lambda=\alpha^2(n+\kappa)-n$, and $\alpha$, $\beta$, and $\kappa$ are constants and $\alpha$ and $\beta$ are nonnegative.

\subsection{Unscented Kalman Filter} \label{ukfp}

Assume the initial estimated mean and the initial estimation error covariance are $\boldsymbol{m}_0$ and $\boldsymbol{P}_0$, 
UKF can be performed in a prediction step and an update step, as in Algorithms \ref{algoKF1} and \ref{algoKF2}.

\begin{algorithm}[!h]
	\caption{UKF Algorithm: Prediction Step}\label{algoKF1}
	\begin{algorithmic}[1]
\State \textbf{calculate} sigma points 
\begin{align} \label{sigmapoint}
&\boldsymbol{\mathcal{X}}_{k-1}=\big [\underbrace{\boldsymbol{m}_{k-1} \cdots \boldsymbol{m}_{k-1}}_{2n+1} \big] \notag \\
&\qquad\qquad + \eta \bmat{\boldsymbol{0}_{n,1} \quad \sqrt{\boldsymbol{P}_{k-1}} \quad -\sqrt{\boldsymbol{P}_{k-1}}\,}.
\end{align}
\State \textbf{evaluate} the sigma points with the dynamic model function
\begin{equation}
\hat{\boldsymbol{\mathcal{X}}}_k=\boldsymbol{f}(\boldsymbol{\mathcal{X}}_{k-1}).
\end{equation}
\State \textbf{estimate} the predicted state mean
\begin{equation}
\boldsymbol{m}_k^-=\sum\limits_{i=0}^{2n}\boldsymbol{W}_m^{(i)}\,\hat{\boldsymbol{\mathcal{X}}}_{i,k}.
\end{equation}
\State \textbf{estimate} the predicted error covariance
\begin{align}
\hspace*{-0.52cm}\boldsymbol{P}_k^-=\sum\limits_{i=0}^{2n}\boldsymbol{W}_c^{(i)}\,(\hat{\boldsymbol{\mathcal{X}}}_{i,k}-\boldsymbol{m}_k^-)
(\hat{\boldsymbol{\mathcal{X}}}_{i,k}-\boldsymbol{m}_k^-)^\top +\boldsymbol{Q}_{k-1}. \label{step4}
\end{align}
\State \textbf{calculate} the predicted sigma points
\begin{align} \label{sigmapoint1}
\hspace*{-0.46cm}\boldsymbol{\boldsymbol{\mathcal{X}}}_k^-=\big [\underbrace{\boldsymbol{m}_{k}^- \cdots \boldsymbol{m}_{k}^-}_{2n+1} \big] + \eta \bmat{\boldsymbol{0}_{n,1} \quad \sqrt{\boldsymbol{P}_{k}^-} \quad -\sqrt{\boldsymbol{P}_{k}^-}\,}.
\end{align}
\State \textbf{evaluate} the propagated sigma points with measurement function
\begin{equation}
\boldsymbol{\boldsymbol{\mathcal{Y}}}_k^-=\boldsymbol{h}(\boldsymbol{\mathcal{X}}_k^-).
\end{equation}
\State \textbf{estimate} the predicted measurement
\begin{equation}
\boldsymbol{y}_k^-=\sum\limits_{i=0}^{2n}\boldsymbol{W}_m^{(i)}\boldsymbol{\boldsymbol{\mathcal{Y}}}_{i,k}^-.
\end{equation}
\end{algorithmic}
\end{algorithm}

\begin{algorithm}[!h]
	\caption{UKF Algorithm: Update Step}\label{algoKF2}
	\begin{algorithmic}[1]
\State \textbf{estimate} the innovation covariance matrix
\begin{align}
\hspace*{-0.45cm}\boldsymbol{P}_{\tilde{\boldsymbol{y}}_k\tilde{\boldsymbol{y}}_k}=\sum\limits_{i=0}^{2n}\boldsymbol{W}_c^{(i)}\big(
\boldsymbol{\mathcal{Y}}_{i,k}^--\boldsymbol{y}_k^-\big)\big(\boldsymbol{\mathcal{Y}}_{i,k}^--\boldsymbol{y}_k^-\big )^\top + \boldsymbol{R}_k. \label{step1}
\end{align}
\State \textbf{estimate} the cross-covariance matrix
\begin{align}
\boldsymbol{P}_{\boldsymbol{x}_k\boldsymbol{y}_k}=\sum\limits_{i=0}^{2n}\boldsymbol{W}_c^{(i)}\big(
\boldsymbol{\mathcal{X}}_{i,k}^--\boldsymbol{m}_k^-\big)\big(\boldsymbol{\mathcal{Y}}_{i,k}^--\boldsymbol{y}_k^-\big )^\top.
\end{align}
\State \textbf{calculate} the Kalman gain 
\begin{equation}
\boldsymbol{K}_k=\boldsymbol{P}_{\boldsymbol{x}_k \boldsymbol{y}_k} \boldsymbol{P}_{\tilde{\boldsymbol{y}}_k \tilde{\boldsymbol{y}}_k}^{-1}. \label{Kk}
\end{equation}
\State \textbf{estimate} the updated state
\begin{equation}
\boldsymbol{m}_k=\boldsymbol{m}_k^-+\boldsymbol{K}_k\big(\boldsymbol{y}_k-\boldsymbol{y}_k^-\big).
\end{equation}
\State \textbf{estimate} the updated error covariance
\begin{equation}
\boldsymbol{P}_k=\boldsymbol{P}_k^- - \boldsymbol{K}_k \boldsymbol{P}_{\tilde{\boldsymbol{y}}_k \tilde{\boldsymbol{y}}_k}\boldsymbol{K}_k^\top. \label{Pk}
\end{equation}
\end{algorithmic}
\end{algorithm}

\section{Unscented Kalman Filter with Enhanced Numerical Stability} \label{ns}

Here, we propose a UKF-GPS method (see Section \ref{ukf-gps}) and introduce five other approaches to enhance the numerical stability of the classic UKF. 
We also summarize and discuss the advantages and disadvantages of these approaches. 

In Section \ref{ukfp}, the estimation error covariance in Algorithm \ref{algoKF1} should be positive semidefinite, because its square root is required in order to obtain the sigma points, as shown in (\ref{sigmapoint}) and (\ref{sigmapoint1}). However, through propagation the estimation error covariance can lose positive semidefiniteness.

As for why the estimation error covariance can lose positive semidefiniteness for the classic UKF, it has been shown in the Appendix III of
\cite{ukf2} that when $\kappa$, a parameter used for unscented transformation, is negative it is possible to calculate a nonpositive semidefinite estimation error covariance. 
As mentioned in \cite{ukf2}, this problem is not uncommon for methods that approximate higher order moments or probability density distributions, as those described in [\ref{ekf1}],  \cite{lose1}, and \cite{estimation_book}.

In \cite{ukf2} a useful heuristic is proposed as $n+\kappa=3$ which can minimize the moments of the standard Gaussian and the sigma points up to the fourth order. 
From (\ref{sigma}) it is seen that the distance of the sigma point from the mean is proportional to $\eta=\sqrt{n+\kappa}$. 
If the UKF procedure follows the heuristic $n+\kappa=3$, the desired dimensional invariance is achieved by canceling the effect of the system dimension $n$, thus avoiding  
the sampling of nonlocal effects that can lead to significant difficulties in worst cases \cite{sr_ukf1,trans}.
However, for a high dimension system with big $n$, the weight of the center point 
\begin{align}
\m W_c^{(0)}&=\frac{\lambda}{n+\lambda}+(1-\alpha^2+\beta) \notag \\
&=2-\alpha^2+\beta-\frac{n}{3\alpha^2}
\end{align}
can be negative. For a typical selection $\alpha=1,\beta=0$, 
$\m W_c^{(0)}=1-n/3$. When $n>3$, $\m W_c^{(0)}$ will be negative and the calculated covariance may become nonpositive semidefinite.

\subsection{EKF/UKF Toolbox Approach}

In EKF/UKF toolbox \cite{ekfukf}, when $\boldsymbol{P}_{k-1}$ or $\boldsymbol{P}_k^-$ is not positive semidefinite, the function `$\operatorname{schol}$', 
which calculates the lower triangular Cholesky factor of a matrix, can still give an output.
The `$\operatorname{schol}$' algorithm can be summarized as 
\begin{align}
s_{jj}&=\boldsymbol{P}_{jj}-\sum\limits_{k=1}^{j-1}\boldsymbol{S}^2_{jk} \\ 
\boldsymbol{S}_{jj} &= \begin{cases}
	\sqrt{s_{jj}}, &\text{if $\,\,\, s > \epsilon$} \\
	0, &\text{otherwise}
  \end{cases} \\
s_{ij}&=\boldsymbol{P}_{ij}-\sum\limits_{k=1}^{j-1}\boldsymbol{S}_{ik}\boldsymbol{S}_{jk} \\  
\boldsymbol{S}_{ij}&=\begin{cases}
	\frac{s_{ij}}{\boldsymbol{S}_{jj}}, &\text{if $\,\,\, \boldsymbol{S}_{jj} > \epsilon$} \\
	0, &\text{otherwise},
  \end{cases}
\end{align}
where $\boldsymbol{P}$ is the covariance matrix, $\boldsymbol{S}$ is the output of the `$\operatorname{schol}$' function, and 
$\epsilon=2.22\times 10^{-16}$ is the distance from 1.0 to the next largest double-precision number in MATLAB.
If a matrix $\boldsymbol{P}$ is positive semidefinite, `$\operatorname{schol}$' can obtain a $\boldsymbol{S}$ matrix such that $\boldsymbol{P}=\boldsymbol{S}\boldsymbol{S}^\top$. When $\boldsymbol{P}$ is positive semidefinite, the `schol' can still get a matrix $\boldsymbol{S}$ but $\boldsymbol{P}=\boldsymbol{S}\boldsymbol{S}^\top$ cannot be satisfied. 
However, by using this $\boldsymbol{S}$ the sigma points can be calculated and the estimation by UKF can at least continue to proceed. 
This approach for enhancing the numerical stability is called ``\textit{UKF-schol}''.

\subsection{Selection of $\kappa$}

When $\kappa$ is negative it is possible to calculate a nonpositive semidefinite estimation error covariance \cite{ukf2}. 
Therefore, in \cite{sr_ukf1} it is suggested to choose $\kappa \ge 0$ to guarantee the positive semidefiniteness of the the covariance matrix. 
Since the specific value of $\kappa$ is not critical, a good default choice is $\kappa=0$ \cite{sr_ukf1}. 
This approach is named as ``\textit{UKF-$\kappa$}''.

When $\kappa=0$, the distance of the sigma point from the mean is proportional to $\sqrt{n}$.
As $n$ increases, the radius of the sphere that bounds all the sigma points also increases \cite{sr_ukf1}. 
Even though the mean and covariance of the prior distribution are still captured correctly, it does so at the cost of possibly sampling nonlocal effects, which
can lead to significant difficulties if the nonlinearities in question are very severe. 
Therefore, although selecting $\kappa=0$ addresses the numerical instability problem in UKF, it picks up the nonlocal sampling problem.

\subsection{Modified UKF}

In \cite{ukf2} a useful heuristic is proposed as $n+\kappa=3$ which can minimize the moments of the standard Gaussian and the sigma points up to the fourth order. 
This means that for a system with $n>3$, $\kappa$ will be negative. 
In order to avoid obtaining a nonpositive, semidefinite covariance when using a negative $\kappa$, a modified UKF is proposed in \cite{ukf2} for which the predicted error covariance in (\ref{step4}) and the innovation covariance matrix in (\ref{step1}) are evaluated about the projected mean as
\begin{align}
\hspace*{-0.3cm}\boldsymbol{P}_k^-&=\sum\limits_{i=0}^{2n}\boldsymbol{W}_c^{(i)}\,(\hat{\boldsymbol{\mathcal{X}}}_{i,k}-\hat{\boldsymbol{\mathcal{X}}}_{0,k})
(\hat{\boldsymbol{\mathcal{X}}}_{i,k}-\hat{\boldsymbol{\mathcal{X}}}_{0,k})^\top +\boldsymbol{Q}_{k-1} \notag \\
&=\sum\limits_{i=1}^{2n}\boldsymbol{W}_c^{(i)}\,(\hat{\boldsymbol{\mathcal{X}}}_{i,k}-\hat{\boldsymbol{\mathcal{X}}}_{0,k})
(\hat{\boldsymbol{\mathcal{X}}}_{i,k}-\hat{\boldsymbol{\mathcal{X}}}_{0,k})^\top +\boldsymbol{Q}_{k-1}
\end{align}
\begin{align}
\hspace*{-0.35cm}\boldsymbol{P}_{\tilde{\boldsymbol{y}}_k\tilde{\boldsymbol{y}}_k}&=\sum\limits_{i=0}^{2n}\boldsymbol{W}_c^{(i)}\big(
\boldsymbol{\mathcal{Y}}_{i,k}^--\boldsymbol{\mathcal{Y}}_{0,k}^-\big)\big(\boldsymbol{\mathcal{Y}}_{i,k}^--\boldsymbol{\mathcal{Y}}_{0,k}^-\big )^\top + \boldsymbol{R}_k \notag \\
&=\sum\limits_{i=1}^{2n}\boldsymbol{W}_c^{(i)}\big(
\boldsymbol{\mathcal{Y}}_{i,k}^--\boldsymbol{\mathcal{Y}}_{0,k}^-\big)\big(\boldsymbol{\mathcal{Y}}_{i,k}^--\boldsymbol{\mathcal{Y}}_{0,k}^-\big )^\top + \boldsymbol{R}_k.
\end{align}

It is shown in \cite{ukf2} that the modified form ensures positive semidefiniteness, and, in the limit $(n+\kappa)\rightarrow 0$, 
the modified UKF is the same as that of the modified, truncated second-order filter \cite{estimation_book}. This approach is called ``\textit{UKF-modified}''.

\subsection{Adding $\Delta\boldsymbol{Q}$}

In \cite{deltaQ} and \cite{deltaQ1}, an extra positive definite matrix $\Delta Q_k$ is added to the predicted covariance matrix in (\ref{step4}) as a slight modification of the UKF to improve the stability of UKF. It is shown that the estimation error of the UKF is bounded if $\Delta Q_k$ is set properly and the stability of UKF is improved. However, the precision of the estimation can be decreased. This approach is called ``\textit{UKF-$\Delta Q$}''. 
Specifically, the predicted error covariance in (\ref{step4}) becomes
\begin{equation}
\boldsymbol{P}_k^-=\sum\limits_{i=0}^{2n}\boldsymbol{W}_m^{(i)}\,(\hat{\boldsymbol{\mathcal{X}}}_{i,k}-\hat{\boldsymbol{\mathcal{X}}}_{i,0})
(\hat{\boldsymbol{\mathcal{X}}}_{i,k}-\hat{\boldsymbol{\mathcal{X}}}_{i,0})^\top + \hat{\boldsymbol{Q}}_{k-1},
\end{equation}
where $\hat{\boldsymbol{Q}}_{k-1}=\boldsymbol{Q}_{k-1} + \Delta Q_{k-1}$. In \cite{deltaQ} a nonlinear system with linear measurement functions are considered and no method is provided to design the additional $\Delta Q_k$ while in \cite{deltaQ1} a nonlinear system with nonlinear measurement functions are considered and a heuristic method is provided to design $\Delta Q_k$.

\subsection{UKF-GPS} \label{ukf-gps}

If $\boldsymbol{P}_{k-1}$ or $\boldsymbol{P}_k^-$ is nonpositive semidefinite, the UKF-GPS will execute the nearest symmetric positive definite (nearPD) algorithm 
(a R function in `$\operatorname{Matrix}$' package \cite{matrix}), as shown in Algorithm \ref{algoKF3}, by which a symmetric positive semidefinite matrix nearest to $\boldsymbol{P}_{k-1}$ or $\boldsymbol{P}_k^-$ 
in Frobenius norm can be obtained. The input $\boldsymbol{X}_0$ can be $\boldsymbol{P}_{k-1}$ or $\boldsymbol{P}_k^-$ and is converted to the output $\boldsymbol{X}$, which guarantees the positive semidefiniteness and substitutes $\boldsymbol{P}_{k-1}$ or $\boldsymbol{P}_k^-$.

The `$\operatorname{nearPD}$' algorithm adapts the modified alternating projections method in [\ref{spd}] and then adds procedures to force positive definiteness by 
`$\operatorname{posdefify}$' (a R function in `$\operatorname{sfsmisc}$' package) \cite{sfsmisc}, and to guarantee symmetric. 
The modified alternating projections method iteratively projects a matrix onto the set $\mathcal{S}=\{\boldsymbol{Y}=\boldsymbol{Y}^\top \in \mathbb{R}^{n\times n}: \boldsymbol{Y} \ge 0\}$ 
by a modified interation due to Dykstra \cite{dykstra} ($\Delta \boldsymbol{S}$ is Dykstra's correction), which incorporates a judiciously chosen correction to each projection that can be interpreted as a normal vector to the corresponding convex set [\ref{spd}]. 
As is mentioned in [\ref{spd}], general results in \cite{boyle} and \cite{han} show that both $\boldsymbol{X}$ and $\boldsymbol{Y}$ converge to the desired nearest covariance matrix as the number of iterations approach infinity. The rate of convergence of Dykstra's algorithm is linear when the sets are subspaces and the constant depends on the angle between the subspaces \cite{conver}. 
To force positive definiteness, the eigenvalues less than $Eps$ are replaced by a positive value $Eps$.

In Algorithm \ref{algoKF3}, `$\operatorname{eig}$' (eigen decomposition), `$\operatorname{max}$', `$\operatorname{sqrt}$' (square root), `$\operatorname{diag}$', `$\LargerCdot$' (element-wise product), `$\times$' (matrix product), 
and `$./$' (element-wise division) are MATLAB functions; 
$\boldsymbol{V}$ is the matrix of eigenvectors, $\boldsymbol{d}$ is the vector of eigenvalues; 
$\boldsymbol{p}$ is the elements that satisfy $\boldsymbol{d}>\tau_{\textrm{eig}} \, \textrm{max}(\boldsymbol{d})$; 
$[\boldsymbol{V}]_{\boldsymbol{p}}$ is the columns of $\boldsymbol{V}$ that belong to $\boldsymbol{p}$; 
$\boldsymbol{d}_{\boldsymbol{p}}$ is the rows of $\boldsymbol{d}$ that belong to $\boldsymbol{p}$; 
and $||\boldsymbol{A}||$ is the Frobenius norm, the matrix norm of an $m\times n$ matrix $\boldsymbol{A}$ with entry $a_{ij}$ defined as
\begin{equation}
||\boldsymbol{A}||=\sqrt{\sum\limits_{i=1}^{m}\sum\limits_{j=1}^{n}|a_{ij}|^2}.
\end{equation}

\begin{algorithm}[!h]
	\caption{$\operatorname{nearPD}$ Algorithm}\label{algoKF3}
	\begin{algorithmic}[1]
\State \textbf{initialize} 

Let \textnormal{$\Delta \boldsymbol{S}=\boldsymbol{0}_{n,n}$.}
\State \textbf{modified alternating projections}
\begin{align}
\textbf{\textrm{do}} \;\;\;\;\;\;& \notag \\
&\boldsymbol{Y} = \boldsymbol{X} \\
&\boldsymbol{R} = \boldsymbol{Y}-\Delta \boldsymbol{S} \\
&\bmat{\boldsymbol{V},\boldsymbol{d}\,} \leftarrow \operatorname{eig}(\boldsymbol{R}) \\
&\boldsymbol{p} \leftarrow\boldsymbol{d}>\tau_{\textrm{eig}}\,\operatorname{max}(\boldsymbol{d}) \\
&\boldsymbol{X} = [\boldsymbol{V}]_{\boldsymbol{p}} \LargerCdot \big [\underbrace{\boldsymbol{d}_{\boldsymbol{p}} \cdots \boldsymbol{d}_{\boldsymbol{p}}}_{n} \big] \times [\boldsymbol{V}]_{\boldsymbol{p}}^\top  \\
&\Delta \boldsymbol{S}=\boldsymbol{X}-\boldsymbol{R} \\
\textbf{\textrm{while}} \;\;& ||\boldsymbol{Y}-\boldsymbol{X}||/||\boldsymbol{X}||>\tau_{\textrm{conv}} \notag
\end{align}
\State \textbf{guarantee} positive definite
\begin{align}
&\qquad [\boldsymbol{V},\boldsymbol{d}] \leftarrow \operatorname{eig}(\boldsymbol{X}) \\
& \qquad Eps\leftarrow \tau_{\textrm{posd}}\,\operatorname{max}(\boldsymbol{d}) \\
&\qquad\boldsymbol{d}(\boldsymbol{d}<Eps)\leftarrow Eps  \\ 
&\qquad\boldsymbol{diagX} \leftarrow \operatorname{diag}(\boldsymbol{X}) \\
&\qquad\boldsymbol{X} = \boldsymbol{V}\operatorname{diag}(\boldsymbol{d})\boldsymbol{V}^\top \\
&\qquad\boldsymbol{D} = \operatorname{sqrt}\big(\operatorname{max}(Eps,\boldsymbol{diagX})./\operatorname{diag}(\boldsymbol{X})\big ) \\
&\qquad\boldsymbol{X} = \operatorname{diag}(\boldsymbol{D}) \times \boldsymbol{X} \LargerCdot \big [\underbrace{\boldsymbol{D} \cdots \boldsymbol{D}}_{n} \big].
\end{align}
\State \textbf{guarantee} symmetric
\begin{equation}
\boldsymbol{X} = \frac{\boldsymbol{X}+\boldsymbol{X}^\top}{2}.
\end{equation}
\end{algorithmic}
\end{algorithm}

\subsection{SR-UKF} \label{sr ukf}
 
The calculation of the new set of sigma points at the prediction step requires taking a matrix square-root of the covariance matrix $\boldsymbol{P}$ by $\boldsymbol{S}\boldsymbol{S}^\top=\boldsymbol{P}$. 
For UKF, while the square-root of $\boldsymbol{P}$ is an integral part, it is actually still the full covariance $\boldsymbol{P}$ that is recursively updated. 
During the propagation, it is possible that $\boldsymbol{P}$ can lose its positive semidefiniteness. 
By contrast, in the implementation of SR-UKF, $\boldsymbol{S}$ is directly propagated, thus avoiding refactorizing $\boldsymbol{P}$ at each step. 
SR-UKF has been applied to power system DSE in \cite{qi,qi1,qi2}. 

SR-UKF can be implemented by Algorithms \ref{algoKF4} and \ref{algoKF5}. 
The filter is initialized by calculating the matrix square-root of the estimation error covariance once 
via a Cholesky factorization as $\boldsymbol{S}_0=\operatorname{schol}\big(\boldsymbol{P}_0\big)$
where `$\operatorname{schol}$' is a function in EKF/UKF Toolbox that calculates the Cholesky factor of a matrix. 
The propagated and updated Cholesky factor is then used in subsequent iterations to directly form the sigma points. 

Correspondingly, (\ref{4-1})--(\ref{4-2}) in step 4 of Algorithm \ref{algoKF4} replace the estimation error covariance update (\ref{step4}) in Algorithm \ref{algoKF1}; 
(\ref{1-1})--(\ref{1-2}) in Step 1 of Algorithm \ref{algoKF5} replace the innovation covariance update (\ref{step1}) in Algorithm \ref{algoKF2}; 
(\ref{Kk_sr}) replaces (\ref{Kk}) for calculating Kalman gain; and (\ref{U})--(\ref{Sk}) replace (\ref{Pk}) by applying $p$ sequential Cholesky downdates to $\boldsymbol{S}_k^-$ 
where $p$ is the number of outputs.

In Algorithms \ref{algoKF4} and \ref{algoKF5}, the `$\operatorname{qr}$' (orthogonal-triangular decomposition) and `$\operatorname{cholupdate}$' (Rank 1 update to Cholesky factorization) are MATLAB functions; 
`s' denotes the sign of $\boldsymbol{W}_c^{(0)}$ and will be `+' if $\boldsymbol{W}_c^{(0)}>0$ and `-' otherwise. 

We first show why (\ref{4-1})--(\ref{4-2}) is equivalent to (\ref{step4}). 
For the matrix in (\ref{4-1}) which is now denoted by $\boldsymbol{A} \in \mathbb{R}^{3n\times n}$ as
\begin{equation}
\boldsymbol{A}=\bigg[\sqrt{\boldsymbol{W}_c^{(1)}}\,\big(\hat{\boldsymbol{\mathcal{X}}}_{1:2n,k}-\big [\underbrace{\boldsymbol{m}_k^- \cdots \boldsymbol{m}_k^-}_{2n}\big]\big) \; \sqrt{\boldsymbol{Q}_{k-1}}\,\bigg]^\top, 
\end{equation}
a QR decomposition can be performed as
\begin{align}
\boldsymbol{A}&=\tilde{\boldsymbol{Q}}\tilde{\boldsymbol{R}}
	=\bmat{
		\tilde{\boldsymbol{Q}}_1 \; \tilde{\boldsymbol{Q}}_2 } \bmat{
			\tilde{\boldsymbol{R}}_1 \\
			\m 0_{2n\times n}	} 
	=\tilde{\m Q}_1 \tilde{\m R}_1,
\end{align}
where $\tilde{\m Q} \in \mathbb{R}^{3n}$, $\tilde{\m Q}_1 \in \mathbb{R}^{3n\times n}$, and $\tilde{\m Q}_2 \in \mathbb{R}^{3n\times 2n}$ 
are all unitary matrices (for a unitary matrix $\m B$, there is $\m B^\top \m B=\m B \m B^\top = \m I$), $\tilde{\m R}_1 \in \mathbb{R}^{n\times n}$ is an upper triangular matrix, 
$\tilde{\m Q}_1 \tilde{\m R}_1$ is called the thin QR factorization \cite{qr1} or reduced QR factorization \cite{qr2}, and there is 
\begin{equation}
\hspace*{-3.9cm}\tilde{\m R}_1^\top \tilde{\m R}_1 = \tilde{\m R}_1^\top \tilde{\m Q}_1^\top \tilde{\m Q}_1 \tilde{\m R}_1 = {\m A}^\top \m A  \notag
\end{equation}
\begin{equation}
\hspace*{1.0cm}=\sum\limits_{i=1}^{2n}\boldsymbol{W}_c^{(i)}\,(\hat{\boldsymbol{\mathcal{X}}}_{i,k}-\boldsymbol{m}_k^-)
 (\hat{\boldsymbol{\mathcal{X}}}_{i,k}-\boldsymbol{m}_k^-)^\top +\boldsymbol{Q}_{k-1} \label{sk-1}.
\end{equation}

The $\m S_k^-$ in (\ref{sk-}) on the left hand side of the arrow is actually $\tilde{\m R}_1$.  
Then for the $\m S_k^-$ on the left hand side of (\ref{cholupdate}) we have
\begin{equation}
\hspace*{-6cm}(\m S_k^-)^\top \m S_k^- \notag
\end{equation}
\begin{align} \label{sk-2}
\hspace*{-0.2cm}=\begin{cases}
	\tilde{\m R}_1^\top \tilde{\m R}_1 + |\m W_c^{(0)}|(\hat{\boldsymbol{\mathcal{X}}}_{0,k}-\boldsymbol{m}_k^-), &\text{if $\,\,\, \m{W}_c^{(0)} > 0$} \\
	\tilde{\m R}_1^\top \tilde{\m R}_1 - |\m W_c^{(0)}|(\hat{\boldsymbol{\mathcal{X}}}_{0,k}-\boldsymbol{m}_k^-), &\text{otherwise}. 
  \end{cases}
\end{align}
From (\ref{sk-1})-(\ref{sk-2}), it is easy to obtain
\begin{align}
&(\m S_k^-)^\top \m S_k^- \notag \\
=&\sum\limits_{i=0}^{2n}\boldsymbol{W}_c^{(i)}\,(\hat{\boldsymbol{\mathcal{X}}}_{i,k}-\boldsymbol{m}_k^-)
 (\hat{\boldsymbol{\mathcal{X}}}_{i,k}-\boldsymbol{m}_k^-)^\top + \boldsymbol{Q}_{k-1} \notag \\
=&\m P_k^-.
\end{align}
By (\ref{4-2}) we convert the upper triangular matrix to a lower triangular matrix and
for $\m S_k^-$ on the left side of (\ref{4-2}) there is 
\begin{equation}
\m S_k^- (\m S_k^-)^\top=\m P_k^-.
\end{equation}

As for why (\ref{1-1})--(\ref{1-2}) can replace (\ref{step1}), it is similar to why (\ref{4-1})--(\ref{4-2}) is equivalent to (\ref{step4}) 
and thus will not be discussed in detail. The relationship between the $\boldsymbol{S}_{\tilde{\boldsymbol{y}}_k}$ obtained from (\ref{1-1})--(\ref{1-2}) 
and the $\boldsymbol{P}_{\tilde{\boldsymbol{y}}_k\tilde{\boldsymbol{y}}_k}$ in (\ref{step1}) can be written as 
\begin{equation}
\boldsymbol{S}_{\tilde{\boldsymbol{y}}_k}\boldsymbol{S}_{\tilde{\boldsymbol{y}}_k}^\top = \boldsymbol{P}_{\tilde{\boldsymbol{y}}_k\tilde{\boldsymbol{y}}_k}
\end{equation}
and therefore the Kalman gain calculated by (\ref{Kk_sr}) is equivalent to the one in (\ref{Kk}).
Then from (\ref{U})--(\ref{Sk}) we have 
\begin{align}
\m S_k^\top \m S_k&=(\m S_k^-)^\top \m S_k^--(\boldsymbol{K}_k\boldsymbol{S}_{\tilde{\boldsymbol{y}}_k})(\boldsymbol{K}_k\boldsymbol{S}_{\tilde{\boldsymbol{y}}_k})^\top \notag \\
&=\m P_k^--\m K_k \boldsymbol{P}_{\tilde{\boldsymbol{y}}_k\tilde{\boldsymbol{y}}_k} \m K_k^\top
\end{align}
which is implemented by applying $p$ sequential Cholesky downdates to $\boldsymbol{S}_k^-$ 
where $p$ is the number of outputs. Each Cholesky downdates uses one column of $\m U$ as the column vector. 
Thus (\ref{U})--(\ref{Sk}) is equivalent to (\ref{Pk}).

\begin{algorithm}[!h]
	\caption{SR-UKF Algorithm: Prediction Step}\label{algoKF4}
	\begin{algorithmic}[1]
\State \textbf{calculate} sigma points 
\begin{align}
&\boldsymbol{\mathcal{X}}_{k-1}=\big [\underbrace{\boldsymbol{m}_{k-1} \cdots \boldsymbol{m}_{k-1}}_{2n+1}\big] \notag \\
&\qquad \qquad\qquad \quad + \eta \big[\boldsymbol{0}_{n,1} \quad \boldsymbol{S}_{k-1} \quad -\boldsymbol{S}_{k-1} \,\big].
\end{align}
\State \textbf{evaluate} sigma points with the dynamic model function
\begin{equation}
\hat{\boldsymbol{\mathcal{X}}}_{k}=\boldsymbol{f}(\boldsymbol{\mathcal{X}}_{k-1}).
\end{equation}
\State \textbf{estimate} the predicted state mean
\begin{equation}
\boldsymbol{m}_k^-=\sum\limits_{i=0}^{2n}\boldsymbol{W}_m^{(i)}\,\hat{\boldsymbol{\mathcal{X}}}_{i,k}.
\end{equation}
\State \textbf{estimate} the predicted square root of error covariance
\begin{equation}
\hspace*{-6.5cm}[\tilde{\boldsymbol{Q}},\boldsymbol{S}_k^-] \leftarrow \;\notag 
\end{equation}
\begin{equation}
\hspace*{-0.66cm}\operatorname{qr}\Bigg(\bigg[\sqrt{\boldsymbol{W}_c^{(1)}}\,\big(\hat{\boldsymbol{\mathcal{X}}}_{1:2n,k}-\big [\underbrace{\boldsymbol{m}_k^- \cdots \boldsymbol{m}_k^-}_{2n}\big]\big) \; \sqrt{\boldsymbol{Q}_{k-1}}\,\bigg]^\top\Bigg) \label{4-1}
\end{equation}
\begin{equation} \label{sk-}
\hspace*{-6.45cm}\boldsymbol{S}_k^- \leftarrow \bmat{
	\boldsymbol{I}_{n}\\
	\boldsymbol{0}
	}\boldsymbol{S}_k^-
\end{equation}
\begin{equation} \label{cholupdate}
\hspace*{-0.46cm}\boldsymbol{S}_k^- \leftarrow \operatorname{cholupdate}\Big(\boldsymbol{S}_k^-,\sqrt{|\boldsymbol{W}_c^{(0)}|}\;\big(\hat{\boldsymbol{\mathcal{X}}}_{0,k}-\boldsymbol{m}_k^-\big),\textrm{`s'}\Big)
\end{equation}
\begin{equation}
\hspace*{-6.6cm} \m S_k^- \leftarrow (\m S_k^-)^\top. \label{4-2}
\end{equation}
\State \textbf{calculate} predicted sigma points 
\begin{align}
\boldsymbol{\mathcal{X}}_k^-=\big [\underbrace{\boldsymbol{m}_k^- \cdots \boldsymbol{m}_k^-}_{2n+1}\big] + \eta \big[\boldsymbol{0}_{n,1} \quad \boldsymbol{S}_k^- \quad -\boldsymbol{S}_k^- \,\big].
\end{align}
\State \textbf{evaluate} the propagated sigma points with measurement function
\begin{equation}
\boldsymbol{\boldsymbol{\mathcal{Y}}}_{k}^-=\boldsymbol{h}(\boldsymbol{\mathcal{X}}_{k}^-).
\end{equation}
\State \textbf{estimate} the predicted measurement
\begin{equation}
\boldsymbol{y}_k^-=\sum\limits_{i=0}^{2n}\boldsymbol{W}_m^{(i)}\boldsymbol{\boldsymbol{\mathcal{Y}}}_{i,k}^-.
\end{equation}
\end{algorithmic}
\end{algorithm}

\begin{algorithm}[!h]
	\caption{SR-UKF Algorithm: Update Step}\label{algoKF5}
	\begin{algorithmic}[1]
\State \textbf{estimate} the innovation covariance matrix
\begin{equation*}
\hspace*{-6.5cm}[\tilde{\boldsymbol{Q}},\boldsymbol{S}_{\tilde{\boldsymbol{y}}_k}] \leftarrow
\end{equation*}
\begin{equation}
\hspace*{-0.9cm}\operatorname{qr}\bigg(\bigg[\sqrt{\boldsymbol{W}_c^{(1)}}\,\big(\boldsymbol{\boldsymbol{\mathcal{Y}}}_{1:2n,k}^--\big [\underbrace{\boldsymbol{y}_k^- \cdots \boldsymbol{y}_k^-}_{2n}\big] \big)\;\;\sqrt{\boldsymbol{R}_k}\,\bigg]^\top \bigg)  \label{1-1}
\end{equation}
\begin{equation}
\hspace*{-6.55cm}\boldsymbol{S}_{\tilde{\boldsymbol{y}}_k} \leftarrow \bmat{
	\boldsymbol{I}_{p}\\
	\boldsymbol{0}
	}\boldsymbol{S}_{\tilde{\boldsymbol{y}}_k}
\end{equation}
\begin{equation}
\hspace*{-0.59cm}\boldsymbol{S}_{\tilde{\boldsymbol{y}}_k} \leftarrow \operatorname{cholupdate}\Big(\boldsymbol{S}_{\tilde{\boldsymbol{y}}_k},\;
\sqrt{|\boldsymbol{W}_c^{(0)}|}\;\big(\boldsymbol{\mathcal{Y}}_{0,k}^--\boldsymbol{y}_k^- \big),\textrm{`s'}\Big)
\end{equation}
\begin{equation}
\hspace*{-6.65cm} \boldsymbol{S}_{\tilde{\boldsymbol{y}}_k} \leftarrow (\boldsymbol{S}_{\tilde{\boldsymbol{y}}_k})^\top.  \label{1-2}
\end{equation}
\State \textbf{estimate} the cross-covariance matrix
\begin{align}
\boldsymbol{P}_{\boldsymbol{x}_k\boldsymbol{y}_k}=\sum\limits_{i=0}^{2n}\boldsymbol{W}_c^{(i)}\big(
\boldsymbol{\mathcal{X}}_{i,k}^--\boldsymbol{m}_k^-\big)\big(\boldsymbol{\mathcal{Y}}_{i,k}^--\boldsymbol{y}_k^-\big )^\top.
\end{align}
\State \textbf{calculate} the Kalman gain 
\begin{equation}
\boldsymbol{K}_k=\boldsymbol{P}_{\boldsymbol{x}_k \boldsymbol{y}_k} \big(\boldsymbol{S}_{\tilde{\boldsymbol{y}}_k}^\top \big)^{-1} \boldsymbol{S}_{\tilde{\boldsymbol{y}}_k}^{-1}. \label{Kk_sr}
\end{equation}
\State \textbf{estimate} the updated state
\begin{equation}
\boldsymbol{m}_k=\boldsymbol{m}_k^-+\boldsymbol{K}_k\big(\boldsymbol{y}_k-\boldsymbol{y}_k^-\big).
\end{equation}
\State \textbf{estimate} the updated square root of error covariance
\begin{align}
&\boldsymbol{U}=\boldsymbol{K}_k\boldsymbol{S}_{\tilde{\boldsymbol{y}}_k}  \label{U} \\
&\boldsymbol{S}_k=\textrm{cholupdate}\big(\boldsymbol{S}_k^-,\boldsymbol{U},\textrm{`-'} \big). \label{Sk}
\end{align}
\end{algorithmic}
\end{algorithm}

\subsection{Summary and Discussion}

The above-mentioned methods are summarized as follows. 
\begin{enumerate}
\item The UKF-schol approach does not solve the problem of the non-positive semidefiniteness of the estimation error covariance 
but is only able to obtain an inaccurate Cholesky factor when the estimation error covariance is not positive semidefinite.
\item The UKF-$\kappa$ approach guarantees the positive semidefiniteness of the estimation error covariance but discards the useful heuristic $n+\kappa=3$ for $n>3$ and also picks up the nonlocal sampling problem. 
\item UKF-modified can also guarantee the positive semidefiniteness of the estimation error covariance. It is shown that under some conditions it is the same as that of the modified, truncated second-order filter \cite{ukf2}.
\item For UKF-$\Delta Q$ approach, it is hard to select a proper extra positive definite matrix. 
The heuristic proposed in \cite{deltaQ1} does not work for the case with non-positive semidefinite estimation error covariance. 
Also, if the process noise covariance is enlarged too much, the precision may be decreased; if it is not sufficiently enlarged, the estimation error covariance can still be 
non-positive semidefinite. It is more reasonable to find the nearest positive semidefinite matrix, as in UKF-GPS. 
\vspace*{-0.4cm}
\item UKF-GPS converts the estimation error covariance to the nearest positive semidefinite matrix whenever it loses positive semidefinateness. 
However, in some cases in order to guarantee positive semidefiniteness the converted positive semidefinite matrix can be not so close to the original one, and may lead to decrease of precision. 
\item SR-UKF intrinsically guarantees the positive semidefiniteness of the estimation error covariance since the square root of the covariance rather than the covariance itself propagates. 
\item As for the implementation based on the classic UKF, UKF-$\kappa$ and UKF-$\Delta Q$ are easier than the others. 
UKF-schol needs to modify the Cholesky factor algorithm, UKF-modified needs to modify the covariance calculation, and UKF-GPS requires to add the `nearPD' algorithm. 
For SR-UKF, it does require more extensive changes of the Kalman filter procedure. 
\item As for calculation efficiency, SR-UKF can be more efficient than other UKF-based methods, mainly because it makes use of 
powerful linear algebra techniques including the orthogonal-triangular decomposition and Cholesky factor updating.
\end{enumerate}

\section{Power System Dynamic State Estimation}  \label{ps dse}

Here, we discuss how different Kalman filters are applied to dynamic state estimation. 
We apply the generator and measurement model in Section \uppercase\expandafter{\romannumeral 3\relax}.C of \cite{qi}, 
which can be used for multi-machine systems and allows both fourth-order transient generator model and second-order classical generator model.
The terminal voltage phasor and terminal current phasor obtained from PMUs are used as the output measurements.

Let $\mathcal{G}_4$ and $\mathcal{G}_2$ respectively denote the set of generators with fourth-order and second-order model. 
The numbers of generators with fourth-order or second-order model, which are also the cardinality of the sets $\mathcal{G}_4$ and $\mathcal{G}_2$, are $g_4$ and $g_2$, respectively. 
Thus the number of states $n=4\,g_4+2\,g_2$.
For generator $i\in \mathcal{G}_4$, the fast sub-transient dynamics and saturation effects are ignored and the generator model is described by the fourth-order differential equations in local $\textrm{d}$-$\textrm{q}$ reference frame:
\begin{subnumcases}{\label{gen model}}
\dot{\delta_i}=\omega_i-\omega_0 \\
\dot{\omega}_i=\frac{\omega_0}{2H_i}\Big(T_{\textrm{m}i}-T_{\textrm{e}i}-\frac{K_{\textrm{D}i}}{\omega_0}(\omega_i-\omega_0)\Big) \\
\dot{e}'_{\textrm{q}i}=\frac{1}{T'_{\textrm{d0}i}}\Big(E_{\textrm{fd}i}-e'_{\textrm{q}i}-(x_{\textrm{d}i}-x'_{\textrm{d}i})\,i_{\textrm{d}i}\Big) \\
\dot{e}'_{\textrm{d}i}=\frac{1}{T'_{\textrm{q0}i}}\Big(-e'_{\textrm{d}i}+(x_{\textrm{q}i}-x'_{\textrm{q}i})\,i_{\textrm{q}i}\Big),
\end{subnumcases}
where $i$ is the generator serial number.

For generator $i\in \mathcal{G}_2$, the generator model is only described by the first two equations of (\ref{gen model}) and the $e'_{\textrm{q}i}$ and $e'_{\textrm{d}i}$ are kept unchanged. 
The set of generators where PMUs are installed is denoted by $\mathcal{G}_{\textrm{P}}$. For generator $i \in \mathcal{G}_{\textrm{P}}$, $E_{\textrm{t}i}=e_{\textrm{R}i}+je_{\textrm{I}i}$ and $I_{\textrm{t}i}=i_{\textrm{R}i}+ji_{\textrm{I}i}$ can be measured and are used as outputs. $T_{\textrm{m}i}$ and $E_{\textrm{fd}i}$ are used as inputs.

The dynamic model (\ref{gen model}) can be rewritten in a general state space form as
\begin{subnumcases} {\label{n1}}
\dot{\boldsymbol{x}}=\boldsymbol{f}_c(\boldsymbol{x},\boldsymbol{u}) \\
\boldsymbol{y}=\boldsymbol{h}_c(\boldsymbol{x},\boldsymbol{u}),
\end{subnumcases}
where the state vector $\boldsymbol{x}$, input vector $\boldsymbol{u}$, and output vector $\boldsymbol{y}$
are respectively
\begin{subequations}
\begin{align}
\boldsymbol{x} &= \bmat{\boldsymbol{\delta}^\top \quad \boldsymbol{\omega}^\top \quad \boldsymbol{e'_{\textbf{q}}}^\top \quad \boldsymbol{e'_{\textbf{d}}}^\top}^\top \\
\boldsymbol{u} &= \bmat{\boldsymbol{T_\textbf{m}}^\top \quad \boldsymbol{E_{\textbf{fd}}}^\top}^\top \\
\boldsymbol{y} &= \bmat{\boldsymbol{e}_\textbf{R}^\top \quad \boldsymbol{e}_\textbf{I}^\top \quad \boldsymbol{i}_\textbf{R}^\top \quad \boldsymbol{i}_\textbf{I}^\top}^\top.
\end{align}
\end{subequations}

\allowdisplaybreaks
The $i_{\textrm{q}i}$, $i_{\textrm{d}i}$, and $T_{\textrm{e}i}$ in (\ref{gen model}) are actually functions of $\boldsymbol{x}$:
\begin{subequations} \label{temp}
\begin{align}
\Psi_{\textrm{R}i}&=e'_{\textrm{d}i}\sin\delta_i+e'_{\textrm{q}i}\cos\delta_i \\
\Psi_{\textrm{I}i}&=e'_{\textrm{q}i}\sin\delta_i-e'_{\textrm{d}i}\cos\delta_i \\
I_{\textrm{t}i}&=\boldsymbol{\overline{Y}}_i(\boldsymbol{\Psi}_{\textbf{R}}+j \boldsymbol{\Psi}_{\textbf{I}}) \\
i_{\textrm{R}i}&= \operatorname{Re}(I_{\textrm{t}i}) \\
i_{\textrm{I}i}&= \operatorname{Im}(I_{\textrm{t}i}) \\
i_{\textrm{q}i}&=\frac{S_\textrm{B}}{S_{\textrm{N}i}}(i_{\textrm{I}i}\sin\delta_i+i_{\textrm{R}i}\cos\delta_i) \\
i_{\textrm{d}i}&=\frac{S_\textrm{B}}{S_{\textrm{N}i}}(i_{\textrm{R}i}\sin\delta_i-i_{\textrm{I}i}\cos\delta_i) \\
e_{\textrm{q}i}&=e'_{\textrm{q}i}-x'_{\textrm{d}i}i_{\textrm{d}i} \\
e_{\textrm{d}i}&=e'_{\textrm{d}i}+x'_{\textrm{q}i}i_{\textrm{q}i} \\
P_{\textrm{e}i} &= e_{\textrm{q}i}i_{\textrm{q}i}+e_{\textrm{d}i}i_{\textrm{d}i} \\
T_{\textrm{e}i} &= \frac{S_\textrm{B}}{S_{\textrm{N}i}} P_{\textrm{e}i}.
\end{align}
\end{subequations}

In \eqref{temp}, the outputs $i_\textrm{R}$ and $i_\textrm{I}$ are written as functions of $\boldsymbol{x}$.
Similarly, the outputs $e_{\textrm{R}i}$ and $e_{\textrm{I}i}$ can also be written as function of $\boldsymbol{x}$:
\begin{subequations}
\begin{align}
e_{\textrm{R}i}&= e_{\textrm{d}i}\sin\delta_i+e_{\textrm{q}i}\cos\delta_i \\
e_{\textrm{I}i}&=e_{\textrm{q}i}\sin\delta_i-e_{\textrm{d}i}\cos\delta_i.
\end{align}
\end{subequations}

Note that we do not consider the dynamics of $T_\textrm{m}$ and $E_\textrm{fd}$ but assume they are constant and known, since 
the main objective of this paper is to discuss techniques that enhance the numerical stability of UKF. 
The dynamic state estimation with unknown inputs ($T_\textrm{m}$ or $E_\textrm{fd}$) has already been discussed in \cite{ekf4}, \cite{zhou} and 
similar discussion under the framework of this paper will be specially investigated elsewhere.

Similar to \cite{qi} and [\ref{zhou}], the continuous models in (\ref{gen model}) can be discretized into their discrete form as
\begin{subnumcases} {\label{n2}}
\boldsymbol{x}_k = \boldsymbol{f}(\boldsymbol{x}_{k-1},\boldsymbol{u}_{k-1}) \\
\boldsymbol{y}_k = \boldsymbol{h}(\boldsymbol{x}_k,\boldsymbol{u}_k),
\end{subnumcases}
where $k$ denotes the time at $k\Delta_t$ and the state transition functions $\boldsymbol{f}$ can be obtained by the modified Euler method [\ref{kunder}] as
\begin{align}
\tilde{\boldsymbol{x}}_k &= \boldsymbol{x}_{k-1}+\boldsymbol{f}_c(\boldsymbol{x}_{k-1},\boldsymbol{u}_{k-1})\Delta t  \\
\tilde{\boldsymbol{f}} &= \frac{\boldsymbol{f}_c(\tilde{\boldsymbol{x}}_k,\boldsymbol{u}_k) + \boldsymbol{f}_c(\boldsymbol{x}_{k-1},\boldsymbol{u}_{k-1})}{2}\\
\boldsymbol{x}_k &= \boldsymbol{x}_{k-1}+\tilde{\boldsymbol{f}}\Delta t.
\end{align}

The model in (\ref{n2}) can be used to perform power system dynamic state estimation with different Kalman filters. 

\section{Simulation Results} \label{result}

Here, the UKF-GPS and SR-UKF are tested on WSCC 3-machine 9-bus system and NPCC 48-machine 140-bus system, which are extracted from Power System Toolbox (PST) \cite{pst}. 
The EKF and classic UKF comes from EKF/UKF toolbox \cite{ekfukf} and the UKF-GPS and SR-UKF algorithms are implemented based on EKF/UKF toolbox.  
All tests are carried out on a 3.2-GHz Intel(R) Core(TM) i7-4790S based desktop.

\subsection{Settings}

The simulation data is generated as follows.

\begin{enumerate}
\item The simulation data is generated by the model presented in Section \ref{ps dse} and the sampling rate is set to be 120 samples per second.
\item In order to generate dynamic response, a three-phase fault is applied at one bus of the branches with the highest line flows 
and is cleared at the near and remote end after $0.05$ and $0.1$ second. 
We do not consider the fault on lines either bus of which is a generator terminal bus because this can lead to the tripping of a generator. 
\item For each measurement, Gaussian noise with variance $0.01^2$ is added. 
\item The sampling rate of the measurements is set to be 60 frames per second to mimic the PMU sampling rate. 
\item Gaussian process noise is added and the corresponding process noice covariance is set as a diagonal matrix, whose diagonal entries are the square of 10\% of the largest state changes, as in \cite{zhou}.
\item For WSCC system, one PMU is installed at the terminal bus of generator 3, and for NPCC system, 24 PMUs are installed at the terminal bus of generators 1, 2, 3, 4, 6, 9, 
10, 12, 13, 14, 16, 18, 19, 20, 21, 27, 28, 31, 32, 35, 36, 38, 44, and 45; the PMU placements are determined by the method in \cite{qi}, which is based on maximizing the determinant of the empirical observability gramian.
\end{enumerate}

The considered filters are set as follows. 

\begin{enumerate}
\item Dynamic state estimation is performed on the post-contingency system on time period $[0,10\,\textrm{s}]$, which starts from the fault clearing.
\item The initial estimated mean of the system state is set to be the pre-contingency state. 
\item For all methods, $\alpha=1$ and $\beta=0$. For UKF-$\kappa$ method $\kappa=0$ and for all the other methods $\kappa=3-n$.
\item The initial estimation error covariance $\boldsymbol{P}_0$ is set as
\begin{equation}
\renewcommand{\arraystretch}{1.6}
\boldsymbol{P}_0 = \left[ \begin{array}{cccc}
		r_\delta^2 \boldsymbol{I}_g & \boldsymbol{0}_{g,g} & \boldsymbol{0}_{g,g_4} & \boldsymbol{0}_{g,g_4} \\
		\boldsymbol{0}_{g,g} & r_\omega^2 \boldsymbol{I}_{g} & \boldsymbol{0}_{g,g_4} & \boldsymbol{0}_{g,g_4} \\
		\boldsymbol{0}_{g_4,g} & \boldsymbol{0}_{g_4,g} & r_{e'_\textrm{q}}^2 \boldsymbol{I}_{g_4} & \boldsymbol{0}_{g_4,g_4} \\
		\boldsymbol{0}_{g_4,g} & \boldsymbol{0}_{g_4,g} & \boldsymbol{0}_{g_4,g_4} & r_{e'_\textrm{d}}^2  \boldsymbol{I}_{g_4}
\end{array} \right],
\end{equation}
where $r_\delta$ and $r_\omega$ are chosen as $0.5\,\pi /180$ and $10^{-3}\omega_0$, and $r_{e'_\textrm{q}}$ and $r_{e'_\textrm{d}}$ are set to be $10^{-3}$.
\item As mentioned before, the covariance for the process noise is set as a diagonal matrix, whose diagonal entries are the square of 10\% of the largest state changes \cite{zhou}.
\item The covariance for the measurement noise is a diagonal matrix, whose diagonal entries are $0.01^2$, as in \cite{zhou}. 
\item For UKF-$\Delta Q$ method, the additional positive definite matrix $\Delta Q$ is set to be $0.005^2 \m I_n$, as suggested in \cite{deltaQ}.
\item For `$\operatorname{nearPD}$', $\tau_{\textrm{eig}}=10^{-6}$ and $\tau_{\textrm{conv}}=\tau_{\textrm{posd}}=10^{-7}$. 
\end{enumerate}

To quantitatively compare the estimation results, we define the following system state estimation error index
\addtolength{\jot}{1em}
\begin{align}
e_x = \sqrt{\frac{\sum\limits_{i=1}^g\sum\limits_{t=1}^{T_s}\big(x_{i,t}^{\textrm{est}}-x_{i,t}^{\textrm{true}}\big)^2}{g\,T_s}}
\end{align}
where $x$ is a type of states and can be $\delta$, $\omega$, $e'_\textrm{q}$, or $e'_\textrm{d}$;
$x_{i,t}^{\textrm{est}}$ is the estimated state and $x_{i,t}^{\textrm{true}}$ is the corresponding true value
for generator $i$ at time step $t$; $T_s$ is the number of time steps.

\subsection{WSCC 3-Machine System}

\begin{figure}[!t]
\centering
\includegraphics[width=3.4in]{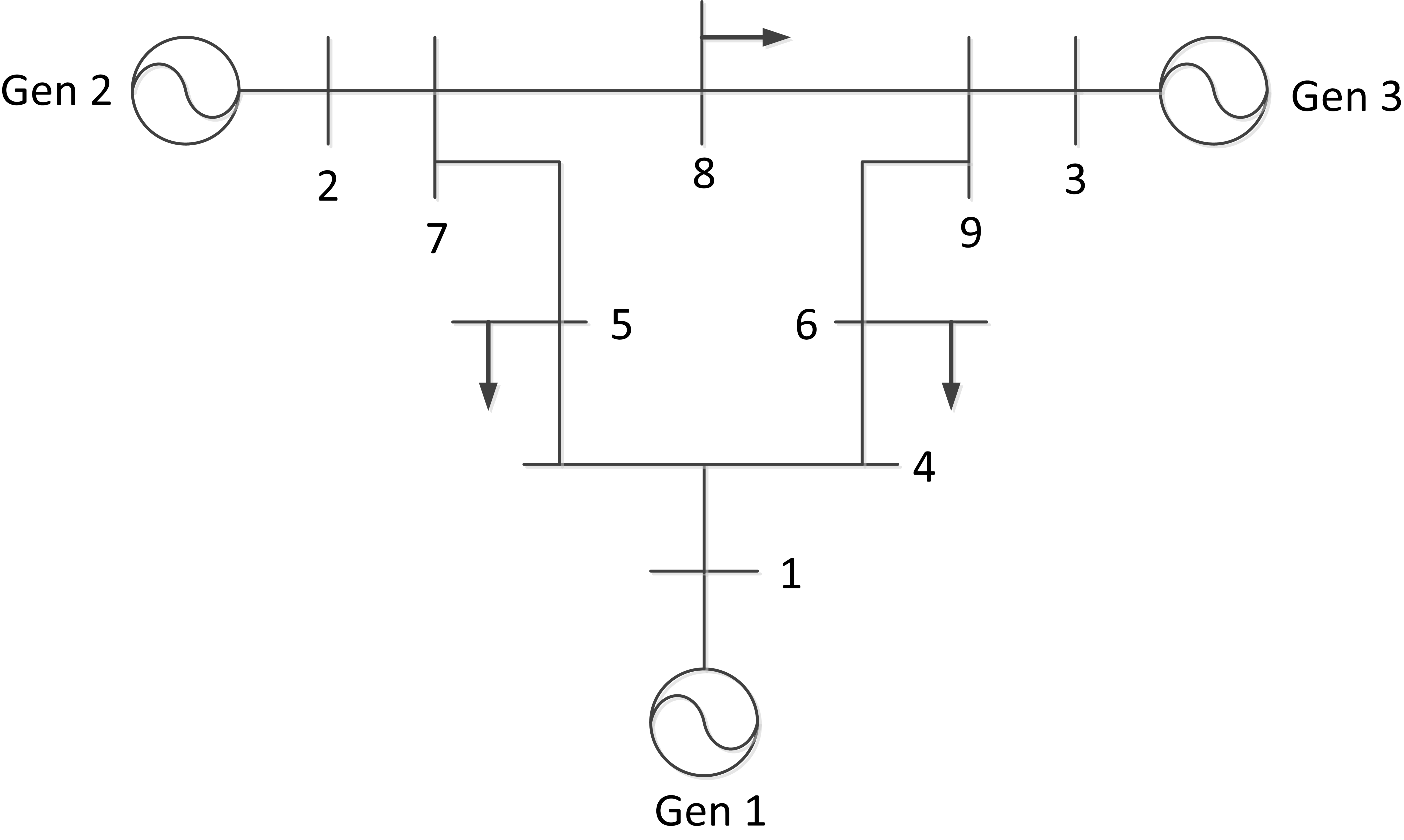}
\caption{WSCC 3-machine 9-bus system.}
\label{wscc3}
\end{figure}

\begin{figure}
\centering
\includegraphics[width=3.49in]{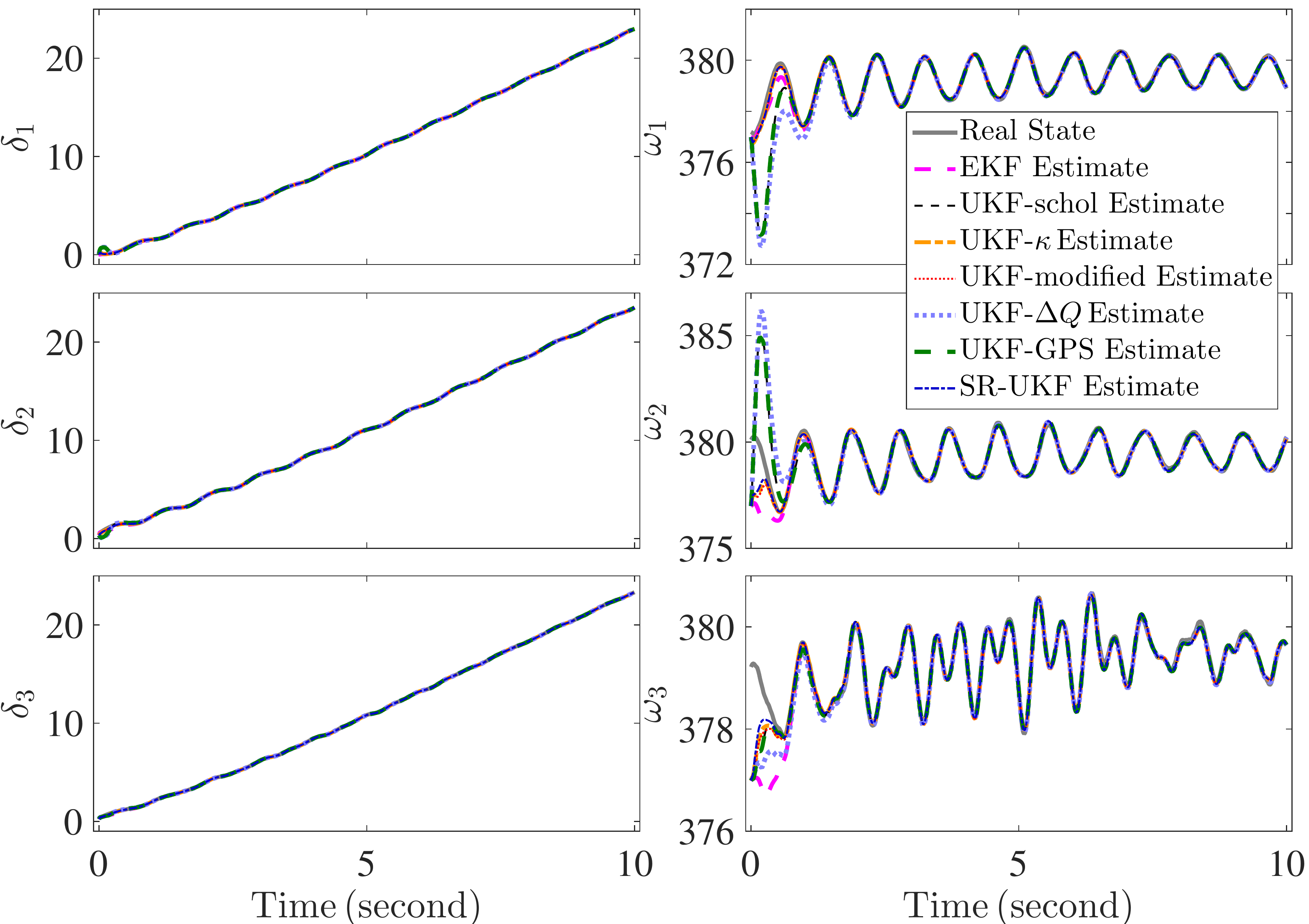}
\caption{Estimated states for WSCC 3-machine system.}
\label{ukf_3}
\end{figure}

Different methods discussed in Section \ref{ns} are tested on the WSCC 3-machine system, as shown in Fig. \ref{wscc3}. 
All generators are assumed to have second-order classical model. 
The estimated state trajectories from different Kalman filters are shown in Fig. \ref{ukf_3}, 
for which a three-phase fault is applied at bus 8 of line $8-9$, the line with the highest line flow. 
For this small system with only six states, there is no obvious numerical stability problem and all methods work well, even though 
for the UKF methods except the UKF-$\kappa$ method there is $\kappa=-3$. 
In this case the estimation error covariance of UKF can keep its positive semidefiniteness during propagation.

There are six branches no bus of which is a generator terminal bus. 
Since the three-phase fault can be applied to any one of the two buses, there are totally twelve possible fault scenarios. 
We perform DSE for each of them and calculate the average values of the system state estimation error index, 
which are listed in Table \ref{error3}. The standard deviations of $\bar{e}_x$ are also listed in the parentheses under $\bar{e}_x$.
It is seen that all methods have small average error and standard deviation and among them SR-UKF has the smallest error and standard deviation.

\begin{table}
\footnotesize
\renewcommand{\arraystretch}{1.3}
\caption{Average Estimation Error for WSCC 3-Machine System}
\label{error3}
\centering
\begin{tabu}{ccc}
\hline
\hline
Filter & $\bar{e}_{\delta}$ & $\bar{e}_{\omega}$ \\
\hline
EKF & \tabincell{c}{0.0371 \\ (0.0167)}  & \tabincell{c}{0.394 \\ (0.0972)}  \\
\hline
UKF-schol & \tabincell{c}{0.0526 \\ (0.0196)} & \tabincell{c}{0.463 \\ (0.159)} \\
\hline
UKF-$\kappa$ & \tabincell{c}{0.0267 \\ (0.0141)} & \tabincell{c}{0.306 \\ (0.102)} \\
\hline
UKF-modified & \tabincell{c}{0.0277 \\ (0.0148)} & \tabincell{c}{0.312 \\ (0.104)} \\
\hline
UKF-$\Delta Q$ & \tabincell{c}{0.0464 \\ (0.0163)} & \tabincell{c}{0.478 \\ (0.148)} \\
\hline
UKF-GPS & \tabincell{c}{0.0526 \\ (0.0196)} & \tabincell{c}{0.463 \\ (0.159)}  \\
\hline
SR-UKF & \tabincell{c}{0.0250 \\ (0.0136)} & \tabincell{c}{0.295 \\ (0.0988)} \\
\hline
\hline
\end{tabu}
\end{table}

\subsection{NPCC 48-Machine System}

\begin{figure}
\centering
\includegraphics[width=3.4in]{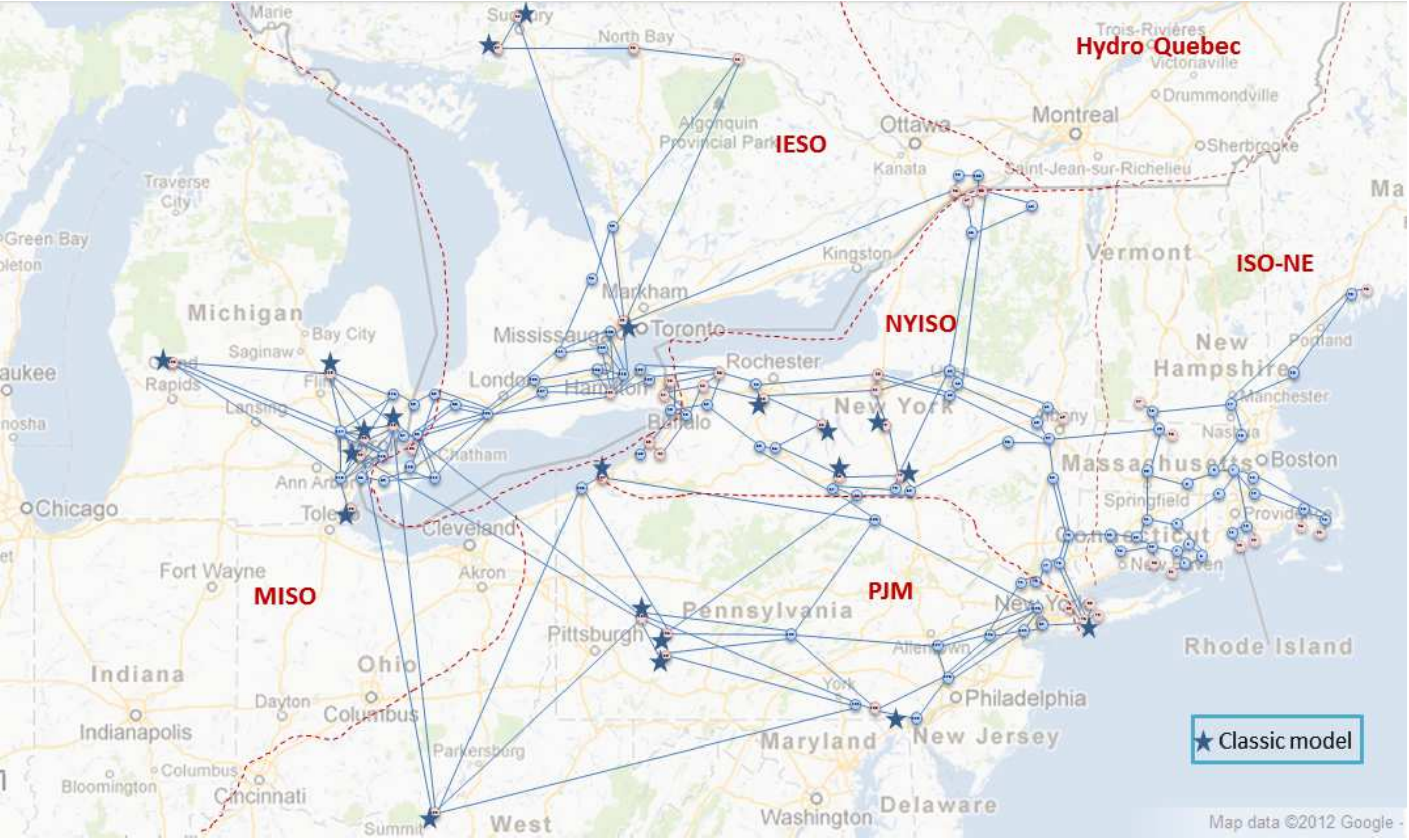}
\caption{Map of the NPCC 48-machine 140-bus system. The stars indicates generators with classical model.}
\label{npcc48}
\end{figure}

As shown in Fig. \ref{npcc48}, the NPCC system [\ref{pst}] represents the northeast region of the EI system. 
Twenty seven generators have fourth-order model and the other twenty one have second-order classical model. 
Thus there are a total of 150 states.

We perform DSE for 50 times and for each of them a three-phase fault is applied at the from bus of one of the 50
branches with highest line flows. For all of the estimations, EKF fails to converge and the classic UKF encounters numerical stability problem 
because the estimation error covariance $\boldsymbol{P}_{k-1}$ or $\boldsymbol{P}_k^-$ loses positive semidefiniteness at some time steps. 
Theoretically, in this case the square root of $\boldsymbol{P}_{k-1}$ or $\boldsymbol{P}_k^-$ cannot be calculated. 
Thus the sigma points in (\ref{sigmapoint}) or (\ref{sigmapoint1}) cannot be obtained and the estimation procedure has to halt. 
Note that both EKF and the classic UKF methods fail due to the infeasibility of the methods themselves rather than other factors such as the settings of the EKF/UKF toolbox or the convergence tolerance: 1) For both methods the EKF/UKF toolbox chooses typical parameters, and using these parameters both methods work well for the smaller WSCC 3-machine system but fail for the bigger NPCC 48-machine system for their poor scalability, which for EKF is because of the loss of nonlinear dynamics in the linearization of the nonlinear transformations and for the classic UKF is due to the above-mentioned numerical stability problem; 2) The estimated states from EKF quickly diverge to values with very large absolute values while the classic UKF cannot continue to perform estimation because of the numerical instability, and thus both methods fail not because of the choice of the convergence tolerance. 

The reason why the estimation error covariance can lose positive semidefiniteness for the classic UKF has been discussed in Section \ref{ns}. 
Here we would like to emphasize that the selection of outputs or the measured values cannot cause the loss of positive semifefiniteness, since we use the same outputs 
and the same settings for simulation data generation and Kalman filters for both WSCC 3-machine system and NPCC 48-machine system and the estimation for WSCC system works very well.  
Also, the measurement configuration cannot be the cause since the numerical stability problem still exists even when all of the generators are installed with PMUs.

\begin{figure*}
	\centering
	\includegraphics[scale=0.37]{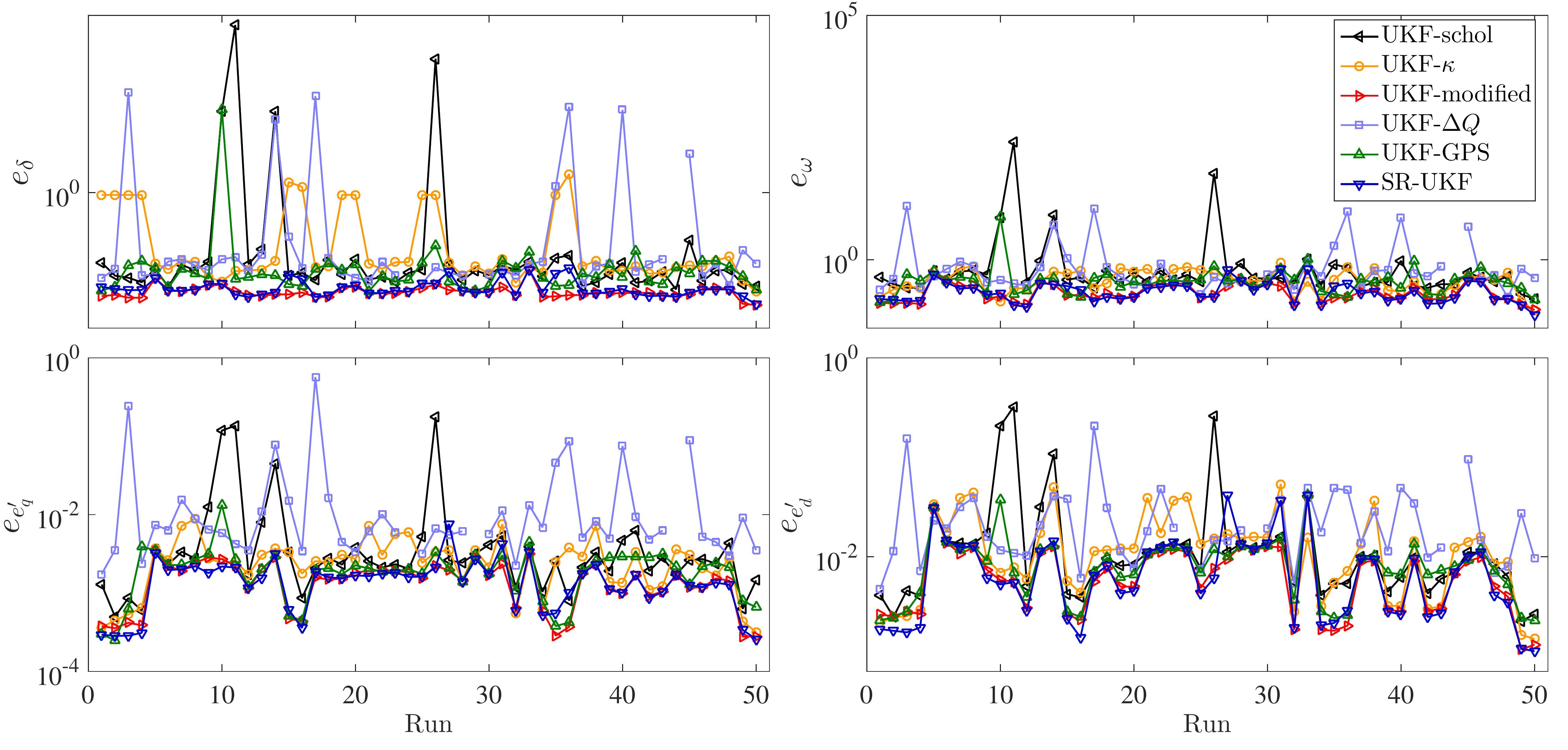}
	\caption{Estimation error index of the states by different methods for NPCC 48-machine system.}
	\label{48m_index}
\end{figure*}

In Fig. \ref{48m_index} we show the estimation error index $e_x$ for each of the fifty estimations. 
We can see that UKF-schol, UKF-$\kappa$, and UKF-$\Delta Q$ do not work well and can have very big estimation errors for several estimations. 
UKF-$\Delta Q$ even diverge for some estimations, for which the estimation error index is too big and thus is not shown.  
UKF-GPS works well for almost all estimations, except for the $10$th estimation in which case it has smaller error for $e'_\textrm{q}$ and $e'_\textrm{d}$ 
than UKF-schol but has similarly big error of $\delta$ and $\omega$. 
By contrast, UKF-modified and SR-UKF both work very well for all estimations due to their enhanced numerical stability and scalability.

\begin{figure}
\centering
\includegraphics[width=3.3in]{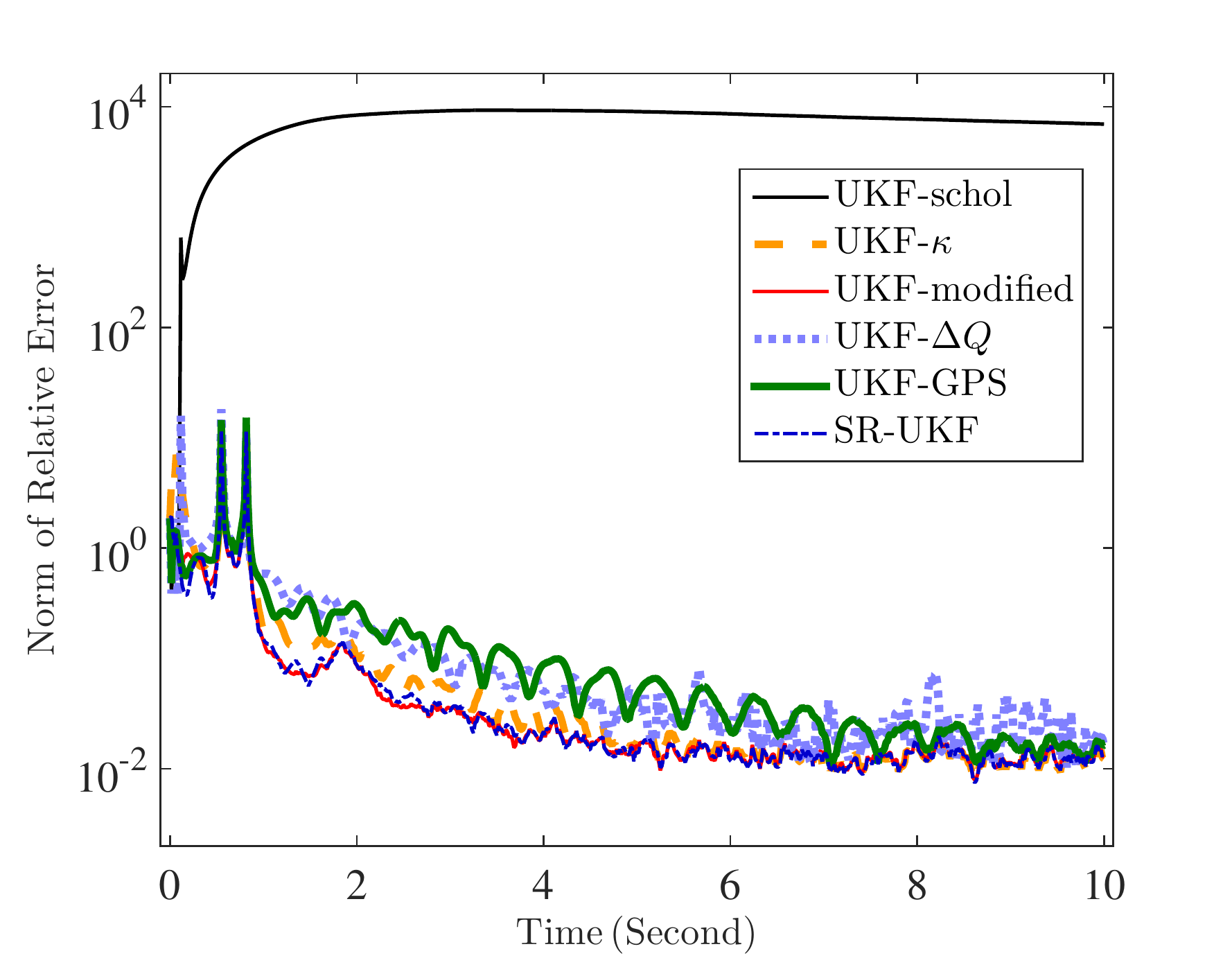}
\caption{Norm of relative error of the states for $11$th estimation.}
\label{est11}
\end{figure}

\begin{figure}
\centering
\includegraphics[width=3.3in]{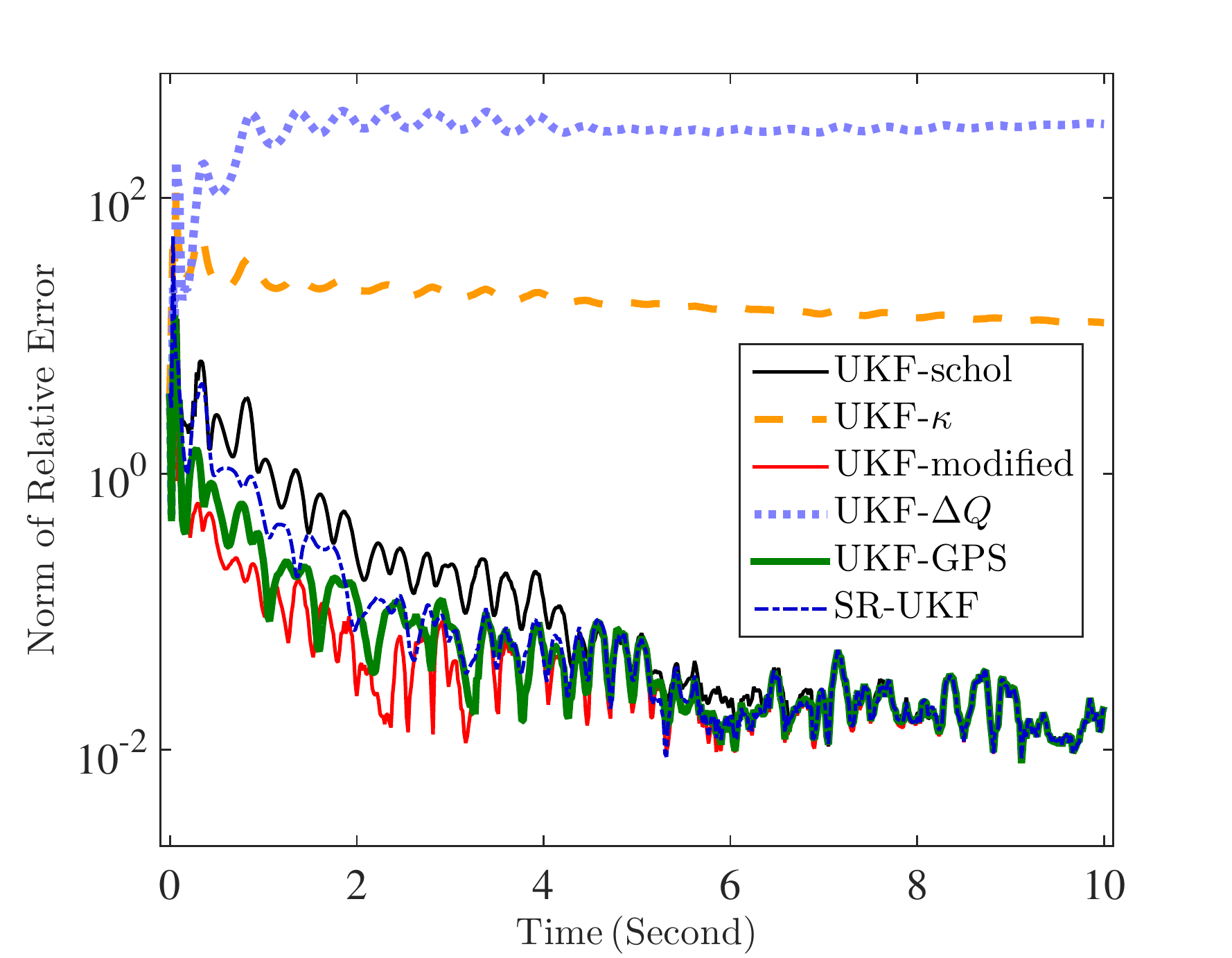}
\caption{Norm of relative error of the states for $36$th estimation.}
\label{est36}
\end{figure}

\begin{figure}
\centering
\includegraphics[width=3.3in]{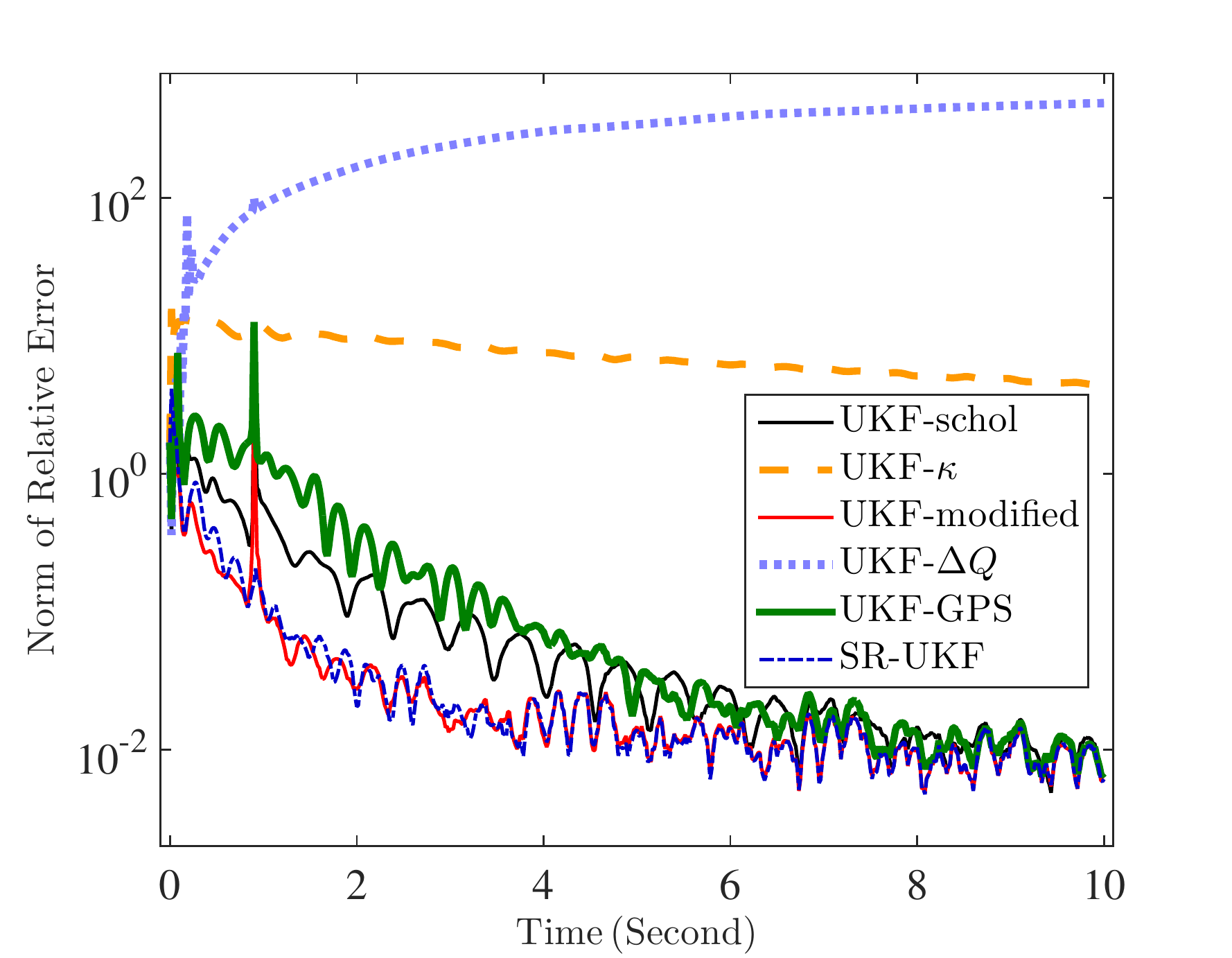}
\caption{Norm of relative error of the states for $3$rd estimation.}
\label{est3}
\end{figure}

\begin{figure}
\centering
\includegraphics[width=3.3in]{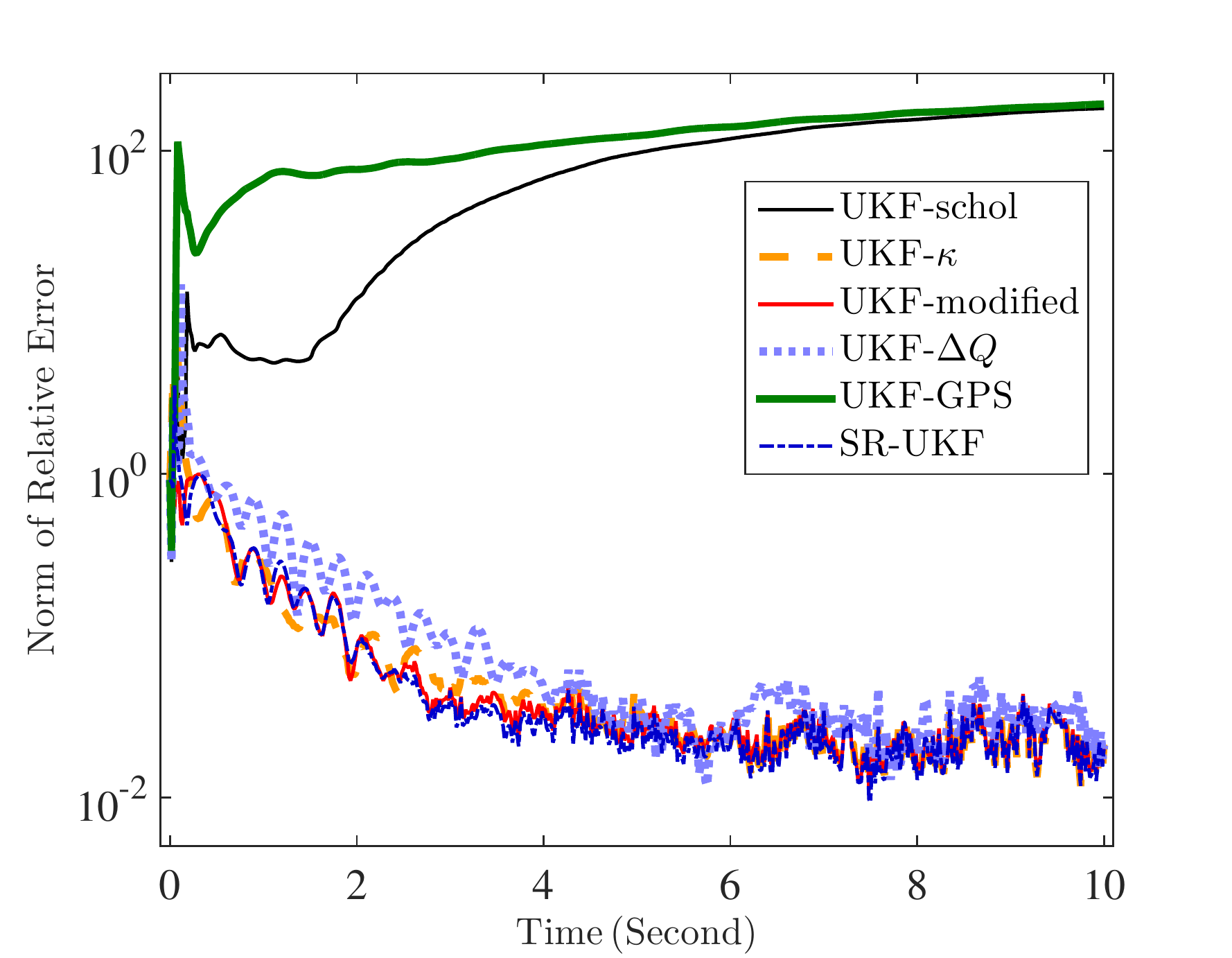}
\caption{Norm of relative error of the states for $10$th estimation.}
\label{est10}
\end{figure}

For the estimation error index of the rotor angle, the UKF-schol, UKF-$\kappa$, UKF-$\Delta Q$, and UKF-GPS get their maximum index among 50 estimations 
on the $11$th, $36$th, $3$rd, $10$th, respectively. 
In Figs. \ref{est11}--\ref{est10}, we show the 2-norm of the relative estimation error of the states $\normof{{(\m x_k -{\m m}_k})/{\m x_k}}_2$ where $\m x_k$ is the real states and $\m m_k$ is the estimated states. From these figures it is seen that the UKF-schol, UKF-$\kappa$, UKF-$\Delta Q$, or UKF-GPS can get poor estimation while the UKF-modified and SR-UKF can always guarantee much better estimation results.

Similar to the WSCC system case, the average values of the estimation error index are also calculated, which are listed in Table \ref{error48}. 
It is seen that the average estimation error index and its standard deviation for UKF-modified and SR-UKF are significantly smaller than the other methods.

\begin{table}
\footnotesize
\renewcommand{\arraystretch}{1.3}
\caption{Average Estimation Error for NPCC 48-Machine System}
\label{error48}
\centering
\begin{tabu}{ccccc}
\hline
\hline
Filter & $\bar{e}_{\delta}$ & $\bar{e}_{\omega}$ & $\bar{e}_{e'_\textrm{q}}$ & $\bar{e}_{e'_\textrm{d}}$ \\
\hline
EKF & --  & --  & -- & -- \\
\hline
UKF-schol & \tabincell{c}{34.018 \\ (191.176)} & \tabincell{c}{6.967 \\ (35.591)} & \tabincell{c}{0.0121 \\ (0.035)} & \tabincell{c}{0.027 \\ (0.063)} \\
\hline
UKF-$\kappa$ & \tabincell{c}{0.295 \\ (0.498)} & \tabincell{c}{0.406 \\ (0.194)} & \tabincell{c}{0.003 \\ (0.002)} & \tabincell{c}{0.015 \\ (0.014)} \\
\hline
UKF-modified & \tabincell{c}{0.0145 \\ (0.004)} & \tabincell{c}{0.232 \\ (0.111)} & \tabincell{c}{0.002 \\ (0.0008)} & \tabincell{c}{0.008 \\ (0.007)} \\
\hline
UKF-$\Delta Q$ & -- & -- & -- & -- \\
\hline
UKF-GPS & \tabincell{c}{0.749 \\ (5.054)} & \tabincell{c}{0.519 \\ (1.027)} & \tabincell{c}{0.002 \\ (0.002)} & \tabincell{c}{0.010 \\ (0.008)} \\
\hline
SR-UKF & \tabincell{c}{0.017 \\ (0.007)} & \tabincell{c}{0.243 \\ (0.132)} & \tabincell{c}{0.002 \\ (0.001)} & \tabincell{c}{0.010 \\ (0.010)} \\
\hline
\hline
\end{tabu}
\end{table}

In the above estimations, we only apply three-phase faults to generate dynamic responses. 
To further validate the proposed approach, we now consider different types of faults, including 
three-phase fault, line to ground fault, line-to-line to ground fault, line-to-line fault, and loss of line.
We perform DSE for 50 times and for each of them a randomly selected type of fault is applied 
at the from bus of one of the 50 branches with highest line flows. 
Similar to the case that only considers three-phase faults, for all estimations EKF fails to converge and classic UKF encounters numerical stability problem. 
The average values of the estimation error index are listed in Table \ref{error48_1}, which shows that 
the UKF-modified and SR-UKF methods have much better performance than the other methods. 
Compared with the case only considering three-phase faults, the estimation error is smaller, possibly because three-phase fault is the most severe fault 
and the corresponding dynamics can be farther away from normal operating conditions.

\begin{table}
\footnotesize
\renewcommand{\arraystretch}{1.3}
\caption{Average Estimation Error for NPCC 48-Machine System under Different Faults}
\label{error48_1}
\centering
\begin{tabu}{ccccc}
\hline
\hline
Filter & $\bar{e}_{\delta}$ & $\bar{e}_{\omega}$ & $\bar{e}_{e'_\textrm{q}}$ & $\bar{e}_{e'_\textrm{d}}$ \\
\hline
EKF & --  & --  & -- & -- \\
\hline
UKF-schol & \tabincell{c}{4.268 \\ (18.236)} & \tabincell{c}{1.274 \\ (3.745)} & \tabincell{c}{0.007 \\ (0.024)} & \tabincell{c}{0.020 \\ (0.060)} \\
\hline
UKF-$\kappa$ & \tabincell{c}{0.174 \\ (0.371)} & \tabincell{c}{0.303 \\ (0.212)} & \tabincell{c}{0.002 \\ (0.002)} & \tabincell{c}{0.010 \\ (0.011)} \\
\hline
UKF-modified & \tabincell{c}{0.013 \\ (0.008)} & \tabincell{c}{0.196 \\ (0.170)} & \tabincell{c}{0.001 \\ (0.001)} & \tabincell{c}{0.008 \\ (0.010)} \\
\hline
UKF-$\Delta Q$ & -- & -- & -- & -- \\
\hline
UKF-GPS & \tabincell{c}{0.586 \\ (3.423)} & \tabincell{c}{0.546 \\ (0.903)} & \tabincell{c}{0.002 \\ (0.002)} & \tabincell{c}{0.010 \\ (0.010)} \\
\hline
SR-UKF & \tabincell{c}{0.013 \\ (0.008)} & \tabincell{c}{0.172 \\ (0.139)} & \tabincell{c}{0.001 \\ (0.001)} & \tabincell{c}{0.006 \\ (0.008)} \\
\hline
\hline
\end{tabu}
\end{table}

As pointed out in \cite{CKF} and \cite{sr_ukf}, EKF, UKF, and SR-UKF all have computational complexity of $\mathcal{O}(n^3)$. 
The average times for performing DSE by different Kalman filters are listed in Table \ref{time}.
Here we list the calculation times for both only considering three-phase fault and randomly choosing different types of faults. 
Note that the time reported here is from MATLAB implementations and is not fully optimized. 
It can be greatly reduced by more efficient, such as C-based, implementations and by further optimization. 
In our implementation the SR-UKF is more efficient than other UKF-based methods, mainly because it makes use of 
powerful linear algebra techniques including the orthogonal-triangular decomposition and Cholesky factor updating. 

It is seen from Table \ref{time} that the additional calculation for `$\operatorname{nearPD}$' is almost negligible 
and the computational complexity of UKF-GPS should also be $\mathcal{O}(n^3)$.
For UKF-GPS, the number of average times that it is requires to execute the `$\operatorname{nearPD}$' algorithm in one estimation 
and the average time steps that need to execute `$\operatorname{nearPD}$' calculation are listed in Table \ref{step}. 
Note that in each time step `$\operatorname{nearPD}$' can be calculated before (\ref{sigmapoint}) or (\ref{sigmapoint1}) in Algorithm \ref{algoKF1} 
and thus the number of times for executing `$\operatorname{nearPD}$' can be greater than the number of time steps involved for `$\operatorname{nearPD}$' calculation.

\begin{table}
\footnotesize
\renewcommand{\arraystretch}{1.3}
\caption{Time for Estimation for NPCC 48-Machine System}
\label{time}
\centering
\begin{tabu}{ccc}
\hline
\hline
\multirow{2}{*}{Filter} & \multicolumn{2}{c}{Time (second)} \\
                         & three-phase fault  &  random fault \\ 
\hline
EKF & 42.615  & 42.661 \\
UKF-schol &  118.580 & 119.246 \\
UKF-$\kappa$ &  118.802 & 119.544. \\
UKF-modified &  118.806 & 119.789 \\
UKF-$\Delta Q$ &  121.188 & 122.274 \\
UKF-GPS & 119.230 & 119.085 \\
SR-UKF & 104.733 & 105.360 \\
\hline
\hline
\end{tabu}
\end{table}

\begin{table}
\footnotesize
\renewcommand{\arraystretch}{1.3}
\caption{Average Times of Executing `nearPD' and Average Number of Time Steps Involved in One Estimation}
\label{step}
\centering
\begin{tabu}{ccc}
\hline
\hline
\multirow{2}{*}{} & \multicolumn{2}{c}{Time (second)} \\
                         & three-phase fault  &  random fault \\ 
\hline
\tabincell{c}{Average times of \\ executing `nearPD'} & 8.88  & 8.50 \\
\hline
\tabincell{c}{Average number of \\ time steps involved}  &  7.34 & 6.90 \\
\hline
\hline
\end{tabu}
\end{table}

\section{Conclusion} \label{conclusion}

In this paper, we introduce and compare six approaches to enhance the numerical stability and further the scalability of the unscented Kalman filter, 
including the proposed UKF-GPS method.
These methods and the extended Kalman Filter are tested on WSCC 3-machine system and NPCC 48-machine system. 
For WSCC system, there is no numerical stability problem for classic UKF, and all methods work well. 
However, for NPCC system, EKF cannot converge and UKF encounters numerical stability problem. 
Among the introduced methods, UKF-schol, UKF-$\kappa$, and UKF-$\Delta Q$ can have big estimation errors for several estimations  
and UKF-$\Delta Q$ even diverge in some cases; UKF-GPS works well for almost all estimations; 
and UKF-modified and SR-UKF work very well for all estimations due to their better numerical stability and scalability. 

Apart from the EKF and UKF that are discussed in this paper, recently some other approaches have also been applied to dynamic state estimation, 
such as the extended particle filter \cite{zhou}, cubature Kalman filter \cite{ckf}, and observers \cite{ckf,obs}. 
EKF, SR-UKF, CKF, and nonlinear observers has been compared for power system DSE under model uncertainty and malicious cyber attacks in \cite{ckf}.
A good comparison of EKF, classic UKF, ensemble Kalman filter, and particle filter is also performed in \cite{zhou1}. 
It would be valuable to more thoroughly compare the approaches discussed in this paper with 
other approaches in order to provide a guideline about how to choose the most suitable approaches for power system DSE.

\begin{IEEEbiography} [{\includegraphics[width=1in,height=1.25in,clip,keepaspectratio]{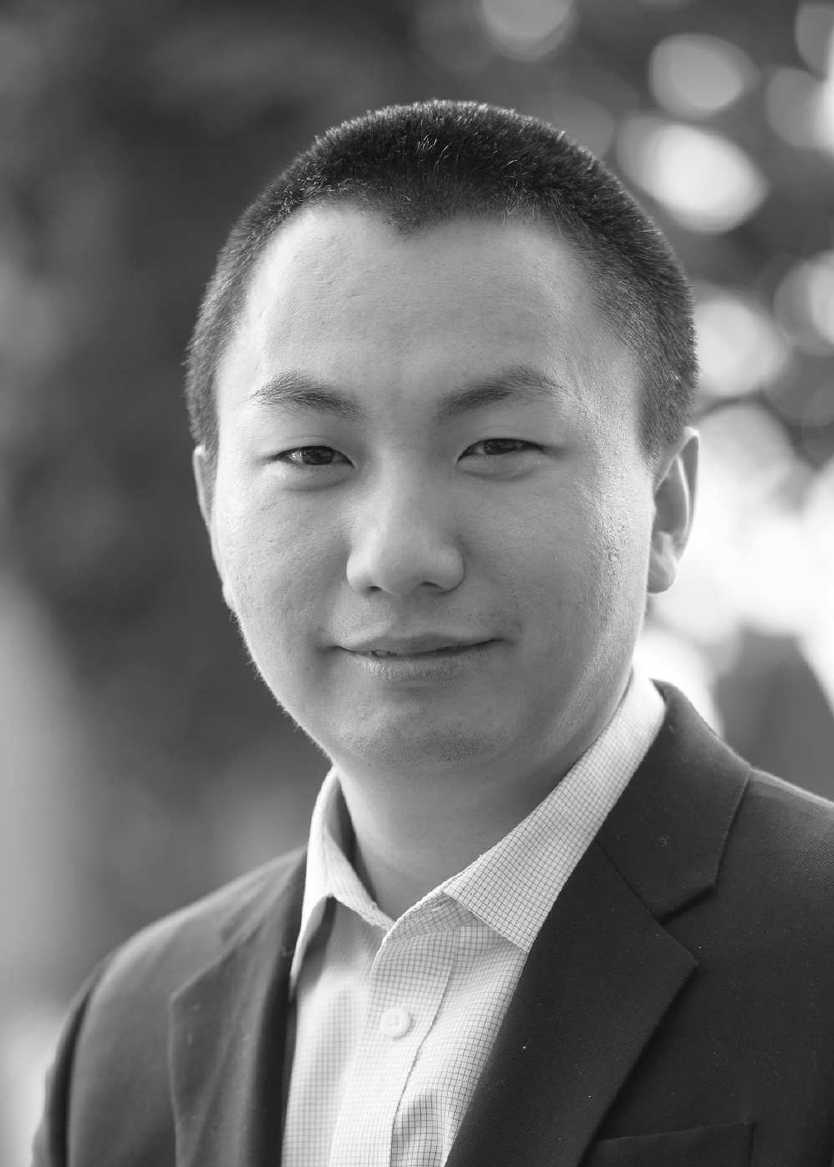}\vfill}]
	{Junjian Qi} (S'12--M'13)
	received the B.E. degree from Shandong University, Jinan, China, in 2008 and the Ph.D. degree Tsinghua University, Beijing, China, in 2013, both in electrical engineering.
	
	In February--August 2012 he was a Visiting Scholar at Iowa State University, Ames, IA, USA. During September 2013--January 2015 he was 
	a Research Associate at Department of Electrical Engineering and Computer Science, University of Tennessee, Knoxville, TN, USA. 
	Currently he is a Postdoctoral Appointee at the Energy Systems Division, Argonne National Laboratory, Argonne, IL, USA. 
	His research interests include cascading blackouts, power system dynamics, state estimation, synchrophasors, and cybersecurity.
\end{IEEEbiography}

\begin{IEEEbiography} [{\includegraphics[width=1in,height=1.25in,clip,keepaspectratio]{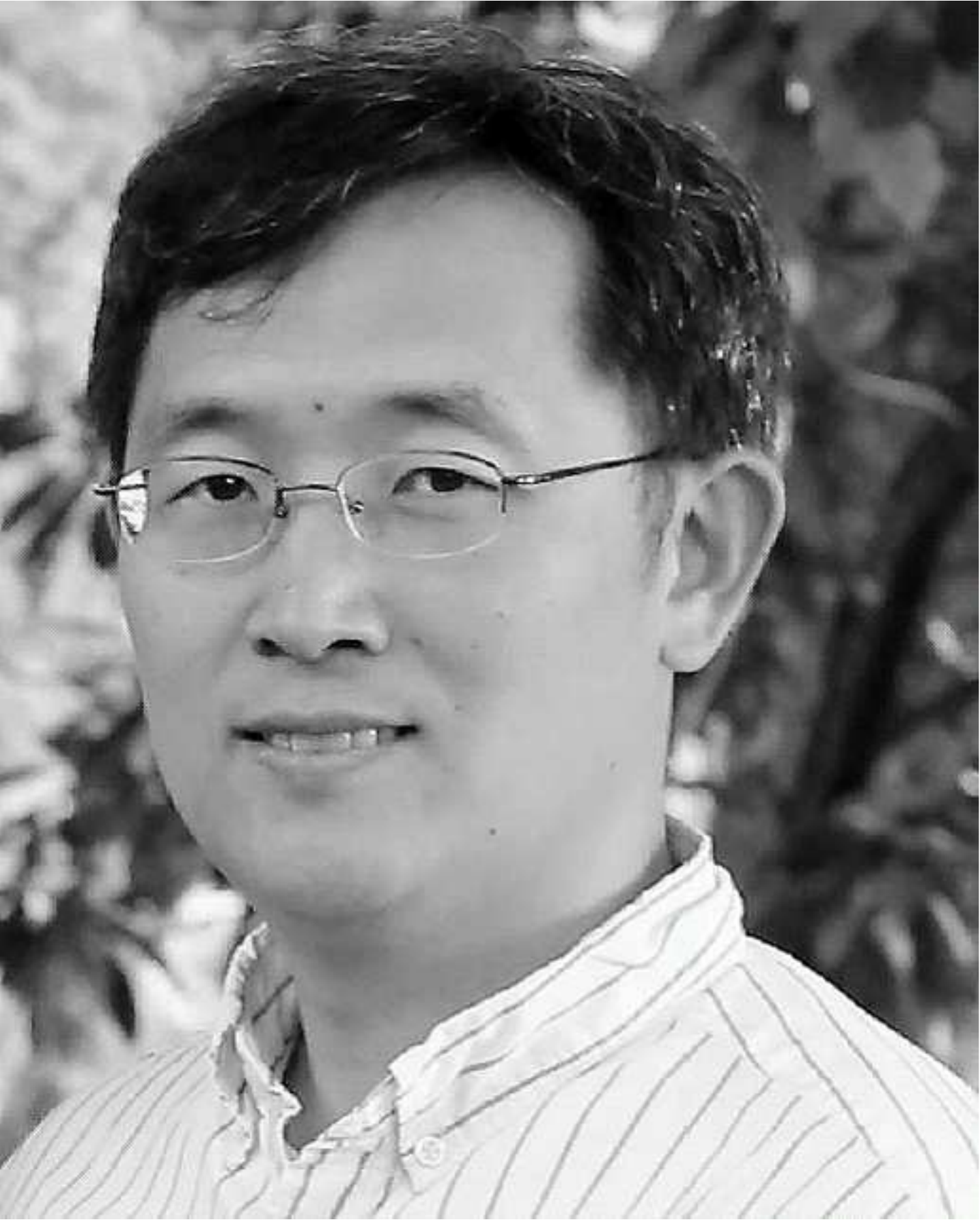}\vfill}]
	{Kai Sun} (M'06--SM'13)
	received the B.S. degree in automation in 1999 and the Ph.D. degree in control science and engineering in 2004 both from Tsinghua University, Beijing, China. 
	
	He is currently an assistant professor at the Department of Electrical Engineering and Computer Science, University of Tennessee in Knoxville. 
	He was a project man-ager in grid operations and planning at the EPRI, Palo Alto, CA from 2007 to 2012. 
	Dr. Sun is an editor of IEEE Transactions on Smart Grid and an associate editor of IET Generation, Transmission and Distribution. 
	His research interests include power system dynamics, stability and control and complex systems.
\end{IEEEbiography}

\begin{IEEEbiography} [{\includegraphics[width=1in,height=1.25in,clip,keepaspectratio]{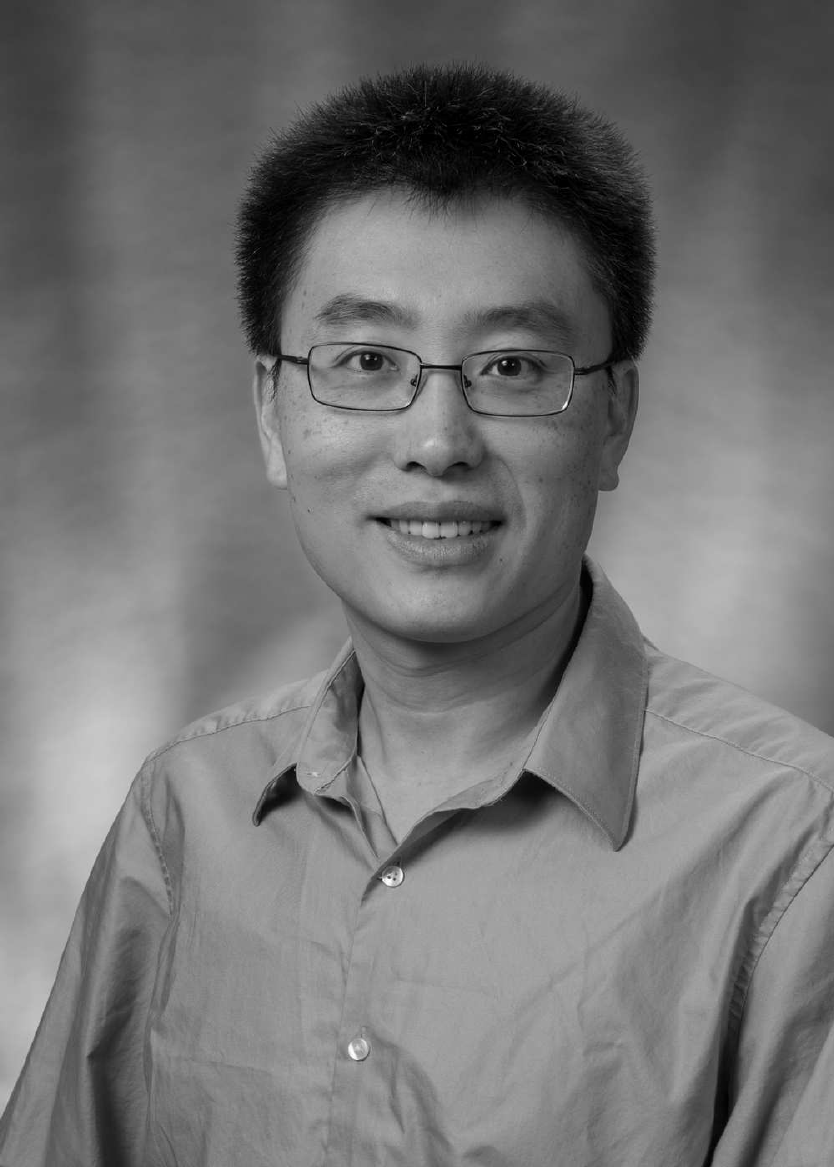}\vfill}]
	{Jianhui Wang} (S'07--SM'12) 
	received the Ph.D. degree in electrical engineering from Illinois Institute of Technology, Chicago, IL, USA, in 2007. 
	
	Presently, he is the Section Lead for Advanced Power Grid Modeling at the Energy Systems Division at Argonne National Laboratory, Argonne, IL, USA.
	Dr. Wang is the secretary of the IEEE Power \& Energy Society (PES) Power System Operations Committee. 
	
	He is an Associate Editor of Journal of Energy Engineering and an editorial board member of Applied Energy. He is also an affiliate professor at Auburn University and an adjunct professor at University of Notre Dame. He has held visiting positions in Europe, Australia, and Hong Kong including a VELUX Visiting Professorship at the Technical University of Denmark (DTU). Dr. Wang is the Editor-in-Chief of the IEEE Transactions on Smart Grid and an IEEE PES Distinguished Lecturer. He is also the recipient of the IEEE PES Power System Operation Committee Prize Paper Award in 2015.
\end{IEEEbiography}

\begin{IEEEbiography} [{\includegraphics[width=1in,height=1.25in,clip,keepaspectratio]{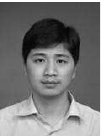}\vfill}]
	{Hui Liu} (M'12)
	received the M.S. degree in 2004 and the Ph.D. degree in 2007 from the School of Electrical Engineering at Guangxi University, China, both in electrical engineering.  
	
	He was a Postdoctoral Fellow at Tsinghua University from 2011 to 2013 and was a staff at Jiangsu University from 2007 to 2016. 
	He visited the Energy Systems Division at Argonne National Laboratory, Argonne, IL, USA, as a visiting scholar from 2014 to 2015. 
	He joined the Department of Electrical Engineering at Guangxi University in 2016, where he is an Associate Professor. 
	His research interests include power system control, electric vehicles, and demand response.
\end{IEEEbiography}


\begin{thebibliography}{11}

\bibitem{Wang1}
B. Wang and K. Sun, ``Power system differential-algebraic equations,'' \emph{arXiv preprint arXiv:1512.05185}, 2015.

\bibitem{se1} \label{se1}
F. C. Schweppe and J. Wildes, ``Power system static-state estimation, Part \uppercase\expandafter{\romannumeral 1\relax}: exact model," \emph{IEEE Trans. Power App. Syst.}, vol. PAS-89, no. 1, pp. 120--125, Jan. 1970.

\bibitem{se2} \label{se2}
A. Abur and A. G\'{o}mez Exp\'{o}sito, \emph{Power System State Estimation: Theory and Implementation}, CRC Press, 2004.

\bibitem{se3} \label{se3}
A. Monticelli, ``Electric power system state estimation," \emph{Proc. IEEE}, vol. 88, no. 2, pp. 262--282, Feb. 2000.

\bibitem{se4} \label{se4}
M. R. Irving, ``Robust state estimation using mixed integer programming," \emph{IEEE Trans. Power Syst.}, vol. 23, no. 3, pp. 1519--1520, Aug. 2008.

\bibitem{se5} \label{se5}
G. He, S. Dong, J. Qi, and Y. Wang, ``Robust state estimator based on maximum normal measurement rate," \emph{IEEE Trans. Power Syst.}, vol. 26, no. 4, pp. 2058--2065, Nov. 2011.

\bibitem{se7} \label{se6}
J. Qi, G. He, S. Mei, and Z. Gu, ``A review of power system robust state estimation,'' \emph{Advanced Technology of Electrical Engineering and Energy}, vol. 30, no. 3, pp. 59--64, Jul. 2011.

\bibitem{se6} \label{se7}
J. Qi, G. He, S. Mei, and F. Liu, ``Power system set membership state estimation," in \emph{Proc. IEEE Power and Energy Soc. Gen. Meeting}, pp. 1--7, San Diego, CA USA, Jul. 2012. 

\bibitem{kf} \label{kf}
R. E. Kalman, ``A new approach to linear filtering and prediction problems," \emph{Trans. ASME J. Basic Eng.}, vol. 82, pp. 34--45, Mar. 1960.

\bibitem{ekf1} \label{ekf1}
A. H. Jazwinski, \emph{Stochastic Processes and Filtering Theory.} San Diego, CA: Academic, 1970.

\bibitem{ekf2} \label{ekf2}
H. W. Sorenson, Ed., \emph{Kalman Filtering: Theory and Application.} Piscataway, NJ: IEEE, 1985.

\bibitem{ekf3} \label{ekf3}
Z. Huang, K. Schneider, and J. Nieplocha, ``Feasibility studies of applying Kalman filter techniques to power system dynamic state estimation," \emph{in Proc. 8th Int. Power Engineering Conf., Singapore}, pp. 376--382, 2007.

\bibitem{ekf4} \label{ekf4}
E. Ghahremani and I. Kamwa, ``Dynamic state estimation in power system by applying the extended Kalman filter with unknown inputs
to phasor measurements," \emph{IEEE Trans. Power Syst.}, vol. 26, no. 4, pp. 2556--2566, Nov. 2011.

\bibitem{CKF} \label{CKF}
I. Arasaratnam and S. Haykin, ``Cubature Kalman filters,'' \emph{IEEE Trans. Autom. Control}, vol. 54, no. 6, pp. 1254--1269, Jun. 2009.

\bibitem{ut} \label{ut}
J. K. Uhlmann, ``Simultaneous map building and localization for real time applications,'' transfer thesis, Univ. Oxford, Oxford, U.K., 1994.

\bibitem{ukf} \label{ukf}
S. J. Julier and J. K.  Uhlmann, ``New extension of the Kalman filter to nonlinear systems," \emph{AeroSense'97}, International Society for Optics and Photonics, pp. 182--193, 1997.

\bibitem{ukf2} \label{ukf2}
S. J. Julier, J. Uhlmann, and H. F. Durrant-Whyte, ``A new method for the nonlinear transformation of means and covariances in filters and estimators,''
\emph{IEEE Trans. Autom. Control}, vol. 45, no. 3, pp. 477--482, Mar. 2000.

\bibitem{ukf1} \label{ukf1}
S. J. Julier and J. K. Uhlmann, ``Unscented filtering and nonlinear estimation," \emph{Proc. IEEE}, vol. 92, no. 3,
pp. 401--422, Mar. 2004.

\bibitem{pwukf1} \label{pwukf1}
E. Ghahremani and I. Kamwa, ``Online state estimation of a synchronous generator using unscented Kalman filter from phasor measurements units,'' \emph{IEEE Trans. Energy Convers.}, vol. 26, no. 4,
pp. 1099--1108, Dec. 2011.

\bibitem{pwukf2} \label{pwukf2}
S. Wang, W. Gao, and A. P. S. Meliopoulos, ``An alternative method for power system dynamic state estimation based on unscented transform," \emph{IEEE Trans. Power Syst.}, vol. 27, no. 2, pp. 942--950, May 2012.

\bibitem{bellman} \label{bellman}
R. E. Bellman, \emph{Adaptive Control Processes}. Princeton, NJ: Princeton Univ. Press, 1961.

\bibitem{ekfukf} \label{ekfukf}
J. Hartikainen, A. Solin, and S. S\"{a}rkk\"{a}, ``Optimal filtering with Kalman filters and smoothers," Dept. of Biomedica Engineering and Computational Sciences, 
Aalto University School of Science, Aug. 2011.

\bibitem{sr_ukf1}
R. Merwe, ``Sigma-point Kalman filters for probabilistic inference in dynamic state-space models,'' PhD diss., Oregon Health \& Science University, 2004.

\bibitem{deltaQ}
K. Xiong, H. Y. Zhang, and C. W. Chan, ``Performance evaluation of UKF-based nonlinear filtering," \emph{Automatica}, vol. 42, no. 2, pp. 261--270, Feb. 2006. 

\bibitem{deltaQ1}
K. Xiong, L. D. Liu, and H. Y. Zhang, ``Modified unscented Kalman filtering and its application in autonomous satellite navigation'', \emph{Aerospace Science and Technology}, vol. 13, no. 4, pp. 238--246, Jul. 2009.

\bibitem{sr_ukf} \label{sr_ukf}
R. Merwe and E. Wan, ``The square-root unscented Kalman filter for state and parameter-estimation," \emph{in Proc. IEEE Int. Conf. Acoustics, Speech, and Signal Processing (ICASSP)}, vol. 6, pp. 3461--3464, 2001. 

\bibitem{lose1}
H. W. Sorenson and A. R. Stubberud, ``Non-linear filtering by approximation of the a posteriori density,'' \emph{Int. J. Contr.}, vol. 8, no. 1, pp. 33--51, Jul. 1968. 

\bibitem{estimation_book}
P. S. Maybeck, \emph{Stochastic Models, Estimation, and Control}, New York: Academic, 1982.

\bibitem{trans}
L. Ghang, B. Hu, A. Li, and F. Qin, ``Transformed unscented Kalman filter,'' \emph{IEEE. Trans. Autom. Control}, vol. 58, no. 1, pp. 252--257, Jan. 2013.

\bibitem{matrix}
D. Bates and M. Maechler, ``Package `Matrix','' Jun. 2015. 

\bibitem{spd} \label{spd}
N. J. Higham, ``Computing the nearest correlation matrix--a problem from finance," \emph{IMA J. Numer. Anal.}, vol. 22, no. 3, pp. 329--343, Jul. 2002.

\bibitem{sfsmisc}
M. Maechler, ``Package `sfsmisc','' Feb. 2015.

\bibitem{dykstra}
R. L. Dykstra, ``An algorithm for restricted least squares regression," \emph{J. Amer. Stat. Assoc.}, vol. 78, no. 384, pp. 837--842, Dec. 1983.

\bibitem{boyle}
J. P. Boyle and R. L. Dykstra, ``A method for finding projections onto the intersection of convex sets in Hilbert spaces,'' \emph{Advances in Order Restricted Inference}, Springer New York, pp. 28--47, 1986.

\bibitem{han}
S. P. Han, ``A successive projection method,'' \emph{Math. Prog.}, vol. 40, no. 1, pp. 1--14, Jan. 1988.

\bibitem{conver}
F. Deutsch and H. Hundal, ``The rate of convergence for the method of alternating projections, II'' \emph{J. Math. Anal. Appl.}, vol. 205, no. 2, pp. 381--405, Jan. 1997. 

\bibitem{qi} \label{qi}
J. Qi, K. Sun, and W. Kang, ``Optimal PMU placement for power system dynamic state estimation by using empirical observability gramian," \emph{IEEE. Trans. Power Syst.}, vol. 30, no. 4, pp. 2041--2054, Jul. 2015.

\bibitem{qi1}
K. Sun, J. Qi, and W. Kang, ``Power system observability and dynamic state estimation for stability monitoring using synchrophasor measurements,'' 
\emph{Control Eng. Pract.}, 2016.

\bibitem{qi2}
J. Qi, K. Sun, and W. Kang, ``Adaptive optimal PMU placement based on empirical observability gramian,'' in \emph{10th IFAC Symposium on Nonlinear Control Systems
(NOLCOS)}, Monterey, CA USA, Aug. 2016.

\bibitem{qr1}
G. H. Golub and  C. F. Van Loan, \emph{Matrix Computations}, JHU Press, 2012.

\bibitem{qr2}
L. N. Trefethen and D. Bau, \emph{Numerical Linear Algebra}, SIAM, 1997.

\bibitem{zhou} \label{zhou}
N. Zhou, D. Meng, and S. Lu, ``Estimation of the dynamic states of synchronous machines using an extended particle filter," \emph{IEEE Trans. Power Syst.}, vol. 28, no. 4, pp. 4152--4161, Nov. 2013.

\bibitem{kunder} \label{kunder}
P. Kunder, \emph{Power System Stability and Control}, New York, NY, USA: McGraw-Hill, 1994.


\bibitem{pst} \label{pst}
J. Chow and G. Rogers, User manual for power system toolbox, Version 3.0, 1991--2008.

\bibitem{ckf}
J. Qi, A. F. Taha, and J. Wang, ``Comparing Kalman filters and observers for dynamic state estimation with model uncertainty and malicious cyber
attacks,'' \emph{arXiv preprint arXiv:1605.01030}, 2016.

\bibitem{obs}
A. F. Taha, J. Qi, J. Wang, and J. H. Panchal, ``Risk mitigation for dynamic state estimation against cyber
attacks and unknown inputs,'' \emph{IEEE Trans. Smart Grid}, to be published.

\bibitem{zhou1}
N. Zhou, D. Meng, Z. Huang, and G. Welch, ``Dynamic state estimation of a synchronous machine using PMU data: A comparative study,'' \emph{IEEE Trans. Smart Grid},
vol. 6, no. 1, pp. 450--460, Jan. 2015.

\end{thebibliography}
\end{document}